\nonstopmode
\documentclass[10pt]{article}

\usepackage{latexsym,amssymb,amsmath}
\usepackage{amsthm, amstext}
\usepackage{array, amsfonts, mathrsfs}
\usepackage{hyperref}
\usepackage{indentfirst}
\usepackage{color}
\usepackage{ifpdf}
\usepackage{url}
\usepackage{mleftright}

\usepackage[T2A]{fontenc}
\usepackage[utf8]{inputenc}
\usepackage[russian,english]{babel}

\usepackage{amsfonts}
\usepackage{amssymb}
\usepackage{amsmath}
\usepackage{pdfpages}
\usepackage{graphicx}
\usepackage{array}
\usepackage{hhline}
\usepackage{doi}
\usepackage{caption}
\usepackage{appendix}

\ifpdf
   \usepackage{pdfpages}
\else
   \RequirePackage[dvips]{graphicx}
\fi

\usepackage{csquotes}

\usepackage[
    backend=biber,
    url=true,
    doi=true,
    eprint=true,
      sorting=nyt,
]{biblatex}

\newcommand{\nocaption}[1]{}

\newcommand{\cchi}{{\bar{\chi}}}

\newcommand{\omegachi}{{\omega}}

\numberwithin{equation}{section}

\newtheorem{myconjecture}{\indent Conjecture}

\newcommand{\myi}{{\mathrm i}}
\newcommand{\myd}{{\mathrm d}}
\newcommand{\mye}{{\mathrm e}}

\begin{document}

\title{Hypothetical
connection of the
theta functions of Dirichlet characters with the real cyclotomic fields}

\author{Yu.\,V.\,Matiyasevich}

\date{}

\maketitle

\begin{quote}\it

We consider a possible approach to the Lindel\"of hypothesis for Dirichlet $L$-functions. It is based on a special form of 
the functional equation for the corresponding
theta functions.
To  estimate $L_\chi(0.5+it)$ we need to solve 
certain systems of linear equations. The entries to the
corresponding matrices  are formed by the summands 
to the series for theta functions and their derivatives. Numerical data suggest that the inverse matrices 
have a deep structure and allow us to state a number of conjectures. In particular,  
it seems that
for a character modulo $q$ certain entries 
to the inverse matrices tend to finite limits
when the sizes of the matrices run over 
arithmetical progressions with step~$2q$. 
Moreover, these limits  belong to the real cyclotomic field~$\mathbb{Q}(\cos(\pi/q))$ (up to a scaling factor of~$\sqrt{q}$). Bibliography 12 items.
  
\end{quote}

\section{Introduction}

One of the most important problems in
Number Theory is about the 
\emph{Riemann zeta function}~$\zeta(s)$.
The celebrated \emph{Riemann Hypothesis}
\cite{Riemann1859}
states that \emph{all non-real zeroes of  $\zeta(s)$
lie on 
the critical line~$\Re(s)=\tfrac{1}{2}$.}
This hypothesis has a lot of important implications
connected with prime numbers.
Today many of such results are proved only conditionally, under the assumption of the Riemann Hypothesis.  

 However, in many cases
it is sufficient to use a weaker assumption,
namely,  the \emph{Lindel\"of Hypothesis} \cite{Lindelof1908}. The latter 
asserts  that
 \emph{for every positive  $\varepsilon$
\begin{equation}\label{Lindelof}
 \zeta(\tfrac{1}{2}+\myi t)=O_\varepsilon(t^\varepsilon)
\end{equation}
as $T$ tends to infinity.}

R. Backlund \cite{Backlund191819} has proved that the Lindelöf Hypothesis is equivalent to the following statement about the zeros of the zeta function: \emph{for every $\varepsilon > 0$, the number of zeros with real part at least $ 1/2 + \varepsilon$ and imaginary part between $T$ and $T + 1$ is $o(\log(T))$
 as $T$ tends to infinity. }
Thus the  Lindelöf Hypothesis can be seen as a weaker form of the Riemann hypothesis.

Both
hypotheses were generalized to
other functions. In particular, it is expected that the
counterparts for these hypotheses hold for all
\emph{Dirichlet $L$-functions}.

A.\,Selberg \cite{Selberg1992} introduced a 
class $S$  of meromorphic
functions. For functions from this class  one can also expects the validity of counterparts
of the Riemann and the Lindel\"of Hypotheses.

Class $S$ is defined by four axioms.
 One of them is 
known as the \emph{functional
equation}. 
It
is considered to be essential  
inevitably necessary 
for the validity of (the counterparts of) the Riemann Hypothesis 
for functions from 
Selberg class --
there are known functions satisfying the other three axioms 
but having zeroes outside the  critical line.

This article originated in an attempt 
to answer the following 
fundamental  question: 
\begin{equation}\label{question}
 \text{\emph{How
the functional equation could be used in the study of
Dirichlet $L$-functions?} }
\end{equation}

Earlier, the author indicated a possible usage of 
the functional equation for proving the Lindel\"of Hypothesis 
for $L$-functions (see \url{http://dx.doi.org/10.13140/RG.2.2.26030.83527}).
In this article we propose a significantly different approach.
It requires solving a particular type of linear systems.
The  entries to the corresponding matrices were suggested by a special
form of the functional equation. 
We  prove that these matrices are not singular.
The inverses to  these matrices are the objects of study in this paper.
Calculations reveal their interesting
properties
and 
 allow us to state a number of conjectures.
 In particular, we conjecture that the limiting 
 values of certain entries to the inverse matrices 
 belong to the real cyclotomic fields.

\section{Functional equations}

Let $q>1$ and $\chi(n)$ be a
\emph{Dirichlet character modulo $q$}.
Corresponding \emph{Dirichlet $L$-function} 
can be defined as
\begin{equation}\label{defL}
  L_\chi(s)=\sum_{n=1}^\infty \chi(n)n^{-s},
\end{equation}
provided that $\Re(s)>1$. We always assume that the character
is not principal, and hence  $L_\chi(s)$ can be extended to an entire function.

There are many ways to state the functional equation. 
Traditionally, this is done via 
\emph{complete $L$-functions}. These functions are defined as follows:
\begin{equation}\label{defxiL}
  \xi_\chi(s)=g_\chi(s)L_\chi(s),
\end{equation}
where
\begin{equation}\label{defgL}
  g_\chi(s)=\left(\tfrac{\pi}{q}\right)^{-\frac{s+\delta}{2}}
  \Gamma\mleft(\tfrac{s+\delta}{2}\mright),
\end{equation}
\begin{equation}\label{defdelta}
\delta=
  \delta(\chi)=
  \begin{cases}
  0,&\text{if }\chi(-1)=1,\\
  1,&\text{if }\chi(-1)=-1.
  \end{cases}
\end{equation}
In terms of the complete $L$-functions, 
the functional equation for a primitive 
character $\chi$ is written as
\begin{equation}\label{feL}
     \xi_\chi(s)=\omegachi\xi_{\bar{\chi}}(1-s),
\end{equation}
where
\begin{equation}\label{defomega}
  \omega(\chi)=\frac{\sum_{k=1}^q\chi(k)
  \mye^\frac{2\pi\myi k}{q}}{\myi^{\delta}\sqrt{q}},
\end{equation}  
and $\bar{\chi}(n)$ is the character conjugate to 
$\chi(n)$,
\begin{equation}\label{defbarchi}
  \bar{\chi}(n)=\overline{\chi(n)}.
\end{equation}

There are many ways to prove the functional equation \eqref{feL}.
After the foundational
 work of
B.\,Riemann \cite{Riemann1859},
  traditionally  this is done via   
 the functional equation for function $\theta_\chi(\tau)$.
 For $\tau$ with positive real part, this function  is defined as follows:
\begin{equation}\label{deftheta}
   \theta_\chi(\tau)=\sum_{n=1}^\infty
  \chi(n)n^\delta \mye^{-\frac{ \pi n  ^2}{q}\tau},
\end{equation}  
where  $\delta$ is the same
as above, namely, is defined by \eqref{defdelta}.
 Function $\theta_\chi(\tau)$ satisfies (see, 
 for example,
  \cite[Ch.\,1, Sect. 4, Lemma\,2]{KarVortrans})  
  \nocite{KarVor}
the functional equation
\begin{equation}\label{thetaFE}
 \theta_\chi({\tau}^{-1})=
\omegachi {\tau^{\delta+\frac{1}{2}}}
  \theta_{\bar{\chi}}(\tau), 
\end{equation}
where $\omega$
 and $\bar\chi$
as above, 
 are defined by 
\eqref{defomega} and \eqref{defbarchi}.

A.\,F.\,Lavrik \cite{Lavrik1991} showed that,
\nocite{Lavrik1991trans}
\emph{vice versa},  
\eqref{thetaFE}   can be deduced from
\eqref{feL}.
Thus  the functional equation~\eqref{thetaFE} can be seen as a faithful form of  
 the functional equation \eqref{feL}
for~$L_\chi(s)$.

The functional equation \eqref{thetaFE} can be restated in a slightly different form. Specifically, this \emph{ functional equality} implies
an infinite series of \emph{numerical equalities}, namely,
\begin{equation}\label{numeqL}
  \left.\frac{\myd^m}{\myd \tau^m}  
 \theta_\chi({\tau}^{-1}) \right\vert_{\tau=1}
 =
  \left.\frac{\myd^m}{\myd \tau^m} 
\omegachi {\tau^{\delta+\frac{1}{2}}}
  \theta_{\bar{\chi}}(\tau)
  \right\vert_{\tau=1},
  \qquad m=0,1,\dots \ . 
\end{equation}
These equalities tell us  that the left- and right-hand sides of
\eqref{thetaFE} have the same derivatives of all orders at $t=1$.
This means that the infinite series of numerical equalities
\eqref{numeqL} implies  the identity \eqref{thetaFE}. Thus, \eqref{numeqL} is yet
another form of the functional equation  \eqref{defxiL} for~$L_\chi(s)$.

Let us introduce some notation. 

For the left-hand side of \eqref{numeqL}, we have:
\begin{eqnarray}
  \left.\frac{\myd^m}{\myd \tau^m} 
 \theta_\chi\mleft({\tau}^{-1}
 \mright)
  \right\vert_{\tau=1}
   &=&
   \frac{\myd^m}{\myd \tau^m}  
  \left.
  \sum_{n=1}^\infty
   \chi(n)n^{\delta}
  \mye^{-\frac{\ \pi n  ^2}{q}\tau^{-1}}
   \right\vert_{\tau=1}
  \\&=&
  \sum_{n=1}^\infty
   \chi(n)n^{\delta}
    \left.
  \frac{\myd^m}{\myd \tau^m}  
  \mye^{-\frac{\ \pi n  ^2}{q}\tau^{-1}}
   \right\vert_{\tau=1}
    \\&=&
  \sum_{n=1}^\infty
   \chi(n)\lambda_{n,m}(\delta,q),
\end{eqnarray}
where 
\begin{eqnarray}\label{deflambda}
  \lambda_{n,m}(d,q)&=&n^{d}
    \left.
  \frac{\myd^m}{\myd \tau^m}  
  \mye^{-\frac{\ \pi n  ^2}{q}\tau^{-1}}
   \right\vert_{\tau=1}\\&=&
 n^{d} 
   \mye^{-\frac{\ \pi n  ^2}{q}}\label{defElambda}
     E_m\mleft(\tfrac{\ \pi n  ^2}{q}\mright)
\end{eqnarray}
for certain polynomials $E_m(x)$.
We can (see Appendix A) give explicit expressions for them:
\begin{eqnarray}\label{defEsum}
  E_{m}(x)&=& \begin{cases}
    1,&\text{ if }m=0,\\
      \sum_{k=1}^m (-1)^{m+k}
  \frac{m!}{k!}\binom{m-1}{k-1}
  x^k,&\text{otherwise.}
  \end{cases}   
\end{eqnarray}

 Similarly,  for the right-hand side of \eqref{numeqL}, we have:
\begin{eqnarray}
 { 
    \left.\frac{\myd^m}{\myd \tau^m} 
       \left( {\tau^{\delta+\frac{1}{2}}}
          \theta_{\bar{\chi}}(\tau)
       \right)
     \right\vert_{\tau=1}\ }
 &=&
 \left.
 \frac{\myd^m}{\myd \tau^m} \left(
  {\tau^{\delta+\frac{1}{2}}}
 \sum_{n=1}^\infty
  \cchi(n)n^{\delta}\mye^{-\frac{\ \pi n  ^2}{q}\tau}
  \right)\right\vert_{\tau=1}
  \\&=&
 \sum_{n=1}^\infty
  \cchi(n)n^{\delta}
  \left.
 \frac{\myd^m}{\myd \tau^m} \left(
  {\tau^{\delta+\frac{1}{2}}}\mye^{-\frac{\ \pi n  ^2}{q}\tau}
  \right)\right\vert_{\tau=1}
    \\&=&
  \sum_{n=1}^\infty
   \cchi(n)\mu_{n,m}(\delta,q),
\end{eqnarray}
where  
\begin{eqnarray}\label{defmu}
  \mu_{n,m}(d,q)&=&n^{d}
  \left.
 \frac{\myd^m}{\myd \tau^m} \left(
  {\tau^{d+\frac{1}{2}}}\mye^{-\frac{\ \pi n  ^2}{q}\tau}
  \right)\right\vert_{\tau=1}\\&=&\label{defFmu}
 n^{d}
    \mye^{-\frac{\ \pi n  ^2}{q}}
      F_{d,m}\mleft(\tfrac{\ \pi n  ^2}{q}\mright)
\end{eqnarray}
for certain polynomials $F_{d,m}(x)$.
For these polynomials, we can (see Appendix B) also give explicit
expressions:
\begin{eqnarray}
  F_{d,m}(x)&=& \label{defFsum2}
    \sum_{k=0}^m(-1)^{k}\binom{m}{k}
 \left(d+\tfrac{1}{2}\right)^{\underline{m-k}}
 x^k;
\end{eqnarray}
 here  $n^{\underline{l}}$ is the \emph{falling 
factorial power},
\begin{equation}\label{fall}
  n^{\underline{l}}=n(n-1)\dots(n-l+1)
\end{equation}
(we use notation from  \cite{ConcMath}).

In the above notation, equalities \eqref{numeqL}
can be rewritten as
\begin{equation}\label{sumEFL}
 \sum_{n=1}^\infty\big(
   \chi(n)\lambda_{n,m}(\delta,q)
     -
   \omegachi  
     \cchi(n)\mu_{n,m}(\delta,q)\big)=0.
   \end{equation}
As  explained above, 
this infinite series (for $m=0,1,\dots$) of numerical equalities is also
equivalent to the functional equation 
 \eqref{defxiL} for~$L_\chi(s)$.

\section{An approach to the Lindel\"of Hypothesis}

The Lindel\"of Hypothesis for $L$-functions asserts that  
 that
 \emph{for every positive  $\varepsilon$
\begin{equation}\label{LindelofL}
 L_\chi(\tfrac{1}{2}+\myi t)=O_\varepsilon(t^\varepsilon)
\end{equation}
as $t$ tends to infinity.}
It can be stated via  the complete 
$L$-function as well. Namely, it is well-known that
 for a positive $a$  
\begin{equation}\label{growthgamma}
  \vert\Gamma(a+\myi\,t)\vert=(1+o_a(t))
  \sqrt{2 \pi }
  t^{a-\frac{1}{2}}
  \mye^{-\frac{ \pi   t} {2}}\quad \text{as }
   t\rightarrow+\infty.
\end{equation}
This implies that \eqref{LindelofL} is equivalent to
 \begin{equation}\label{Lindelofxi}
  \xi_\chi(\tfrac{1}{2}+\myi t)=O_\epsilon\mleft(
  t^{\frac{\delta}{2}-\frac{1}{4}+\epsilon}
  \mye^{-\frac{\pi    t}{4}}\mright).
\end{equation}

According to \eqref{feL}, \eqref{defxiL}, and \eqref{defL},
\begin{eqnarray}\label{Lindelofxi2}
  2\xi_\chi(\tfrac{1}{2}+\myi t)&=&
  \xi_\chi(\tfrac{1}{2}+\myi t)+
  \omega\xi_{\bar\chi}(\tfrac{1}{2}-\myi t)
  \\\label{Lindelofxi3}
  &=&g(\tfrac{1}{2}+  \myi t)L_chi(\tfrac{1}{2}+  \myi t)+
  \omega g(\tfrac{1}{2}+  \myi t)L_chi(\tfrac{1}{2}+  \myi t)
  \\\label{Lindelofxi4}
  &=&\sum_{n=1}^\infty\left(\chi(n)g(\tfrac{1}{2}+  \myi t)
  n^{-({1}/{2}- \myi t)}+
\omega\bar\chi(n)g(\tfrac{1}{2}- \myi t)
  n^{-({1}/{2}- \myi t)}\right).
  \end{eqnarray}
Thus, in order to prove the Lindel\"of Hypothesis, 
we need to estimate the absolute value of \eqref{Lindelofxi4}
from above. 

As an answer to the fundamental question \eqref{question},
here 
we suggest the following possible approach to the Lindel\"of Hypothesis. 
The sum  \eqref{Lindelofxi4}  has the structure similar to
the structure of the sum in the left-hand side of \eqref{sumEFL}. Thanks to the functional equation, we know that the 
latter sum vanishes.
 We can try to use this knowledge to estimate \eqref{Lindelofxi4}.

For example, we can  formally write
\begin{multline}\label{Lindelofxi5}
  2\xi_\chi(\tfrac{1}{2}+\myi t)
  =\sum_{n=1}^\infty\left(\chi(n)
  \bigg(g(\tfrac{1}{2}+  \myi t)
  n^{-({1}/{2}+  \myi t)}+\sum_{m=0}^\infty 
  \nu_{\delta,q,m}(t)\lambda_{n,m}(\delta,q)\bigg)\right.+\\
+\left.\omega\bar\chi(n)
  \bigg(g(\tfrac{1}{2}-  \myi t)
  n^{-({1}/{2}-  \myi t)}-\sum_{m=0}^\infty 
  \nu_{\delta,q,m}(t)\mu_{n,m}(\delta,q)\bigg)\right)
  \end{multline}
 with any weights $\nu_{\delta,q,m}(t)$.
 They 
 need to be chosen so that all the summands in 
 \eqref{Lindelofxi5} are sufficiently small (in absolute value).
 
  In this paper 
 we consider one particular choice of $\nu_{\delta,q,m}(t)$. We define
   them 
 depending on an extra integer parameter~$N$,
 \begin{equation}
    \nu_{\delta,q,m}(t)= \nu_{N,\delta,q,m}(t).
 \end{equation}
For $m\ge 2N$ we put 
 \begin{equation}
   \nu_{N,\delta,q,m}(t)=0
 \end{equation}
  and the
   initial (for $m=0,\dots, 2N-1$)  weights
   should satisfy the equalities  
 \begin{equation}\label{eql}
 g(\tfrac{1}{2}+  \myi t)
  n^{-({1}/{2}+  \myi t)}+\sum_{m=0}^{2N-1} 
  \nu_{N,\delta,q,m}(t)\lambda_{n,m}(\delta,q)=0,
   \end{equation}
 \begin{equation}\label{eqm}
 g(\tfrac{1}{2}-  \myi t)
  n^{-({1}/{2}-  \myi t)}-\sum_{m=0}^{2N-1}
  \nu_{N,\delta,q,m}(t)\mu_{n,m}(\delta,q)=0
 \end{equation} 
   for $n=1,\dots,N$.
  In such a case
  \begin{eqnarray}\label{Lindelofxi6}
  \lefteqn{2\xi_\chi(\tfrac{1}{2}+\myi t)
  =}
  \\&=&\label{Lindelofxi7}
  \sum_{n=N+1}^\infty\left(\chi(n)
  \bigg(g(\tfrac{1}{2}+  \myi t)
  n^{-({1}/{2}+  \myi t)}+\sum_{m=0}^{2N-1} 
  \nu_{N,\delta,q,m}(t)\lambda_{n,m}(\delta,q)\bigg)\right.+
  \\&&\qquad +
  \left.\omega\bar\chi(n)
  \bigg(g(\tfrac{1}{2}-  \myi t)
  n^{-({1}/{2}-  \myi t)}-\sum_{m=0}^{2N-1}
  \nu_{N,\delta,q,m}(t)\mu_{n,m}(\delta,q)\bigg)\right)
\\&=& \label{Lindelofxi8}
  \sum_{n=N+1}^\infty\bigg(\chi(n)
  g(\tfrac{1}{2}+  \myi t)
  n^{-({1}/{2}+  \myi t)}+
\omega\bar\chi(n)
  g(\tfrac{1}{2}-  \myi t)
  n^{-({1}/{2}-  \myi t)}\bigg)
+\\&&\qquad + \label{Lindelofxi9}
\sum_{m=0}^{2N-1}\nu_{N,\delta,q,m}(t)
  \sum_{n=N+1}^\infty\bigg(\chi(n)
  \lambda_{n,m}(\delta,q)-
\omega\bar\chi(n)
    \mu_{n,m}(\delta,q)\bigg).
  \end{eqnarray}
  
 The sum \eqref{Lindelofxi8} is a tail of 
 the convergent sum  \eqref{Lindelofxi4};
 therefore, \eqref{Lindelofxi8} can be made arbitrary small by selecting $N$ sufficiently large.

 Similarly, the inner sum in \eqref{Lindelofxi9} 
 is a tail of 
 the convergent sum  \eqref{sumEFL};
 therefore, for a fixed $m$ each such inner sum
 vanishes as $N$ tends to $\infty$.

 Now our success of estimating $\xi_\chi(\tfrac{1}{2}+\myi t)$ 
 depends on the rate of growth of numbers 
 $\nu_{N,\delta,q,m}$. These numbers were defined via 
 the system of linear equations \eqref{eql}--\eqref{eqm}.
  This system can be written 
 in
 matrix notation as follows.

Let
 \begin{equation}\label{defB}
   B_N(d,q)=\left[
   \begin{array}{ccc}
      \lambda_{N,0}(d,q)
           &\dots&
      \lambda_{N,2N-1}(d,q)
     \\\vdots&\ddots&\vdots\\
         \lambda_{1,0}(d,q)
           &\dots&
      \lambda_{1,2N-1}(d,q)
      \\  \mu_{1,0}(d,q)
           &\dots&
      \mu_{1,2N-1}(d,q)
     \\\vdots&\ddots&\vdots\\
        \mu_{N,0}(d,q)
           &\dots&
      \mu_{N,2N-1}(d,q)
   \end{array}
   \right].
 \end{equation}
(Matrix $B_N({d,q})$ is a submatrix of matrix   $B_{N+1}({d,q})$, and thus can be viewed 
 as a submatrix of the infinite matrix  $B_\infty({d,q})$.)
 In this notation, system \eqref{eql}--\eqref{eqm}
 can be written as
 \begin{equation}
  B_N(d,q)   
    \left[\begin{array}{c}
     \nu_{N,d,q,0}(t)\\
     \vdots\\\vdots\\[1mm]\vdots\\\vdots\\ 
      \nu_{N,d,q,2N-1}(t)
   \end{array}     \right]
   =
   \left[\begin{array}{c}
     \phantom{-}g(\tfrac{1}{2}+  \myi t)
  N^{-({1}/{2}+  \myi t)}\\ \vspace{-2mm}
     \vdots\\\\
       \phantom{-}g(\tfrac{1}{2}+  \myi t)
  1^{-({1}/{2}+  \myi t)}\\[1mm]
       -g(\tfrac{1}{2}-  \myi t)
  1^{-({1}/{2}- \myi t)}\\  \vspace{-2mm} \vdots\\
        \\
       -g(\tfrac{1}{2}-  \myi t)
  N^{-({1}/{2}-  \myi t)}
   \end{array}     \right].
 \end{equation}

 For a rational $d$, matrices $B_N({d,q})$ are not singular (see  Appendix C).
 Therefore, for $d=0$ and for $d=1$ 
 numbers $\nu_{N,d,q,m}$ are well-defined:
 \begin{equation}
     \left[\begin{array}{c}
     \nu_{N,d,q,0}(t)\\
     \vdots\\\vdots\\[1mm]\vdots\\\vdots\\ 
      \nu_{N,d,q,2N-1}(t)
   \end{array}     \right]
   =B_{N}(d,q)^{-1}
   \left[\begin{array}{c}
     \phantom{-} g(\tfrac{1}{2}+  \myi t)
  N^{-({1}/{2}+  \myi t)}\\ \vspace{-2mm}
     \vdots\\\\
      \phantom{-} g(\tfrac{1}{2}+  \myi t)
  1^{-({1}/{2}+  \myi t)}\\[1mm]
       -g(\tfrac{1}{2}-  \myi t)
  1^{-({1}/{2}- \myi t)}\\  \vspace{-2mm} \vdots\\\\
       -g(\tfrac{1}{2}-  \myi t)
  N^{-({1}/{2}-  \myi t)}
   \end{array}     \right].
 \end{equation}

 The inverse matrices $B_N^{-1}({d,q})$
 for $d=0$ and for $d=1$ 
 are the main object of investigation in this paper.
 They have deep internal structure 
 and here we study them \emph{per se};
 above described 
 possible application to the Lindel\"of Hypothesis
 will be considered in another paper.
 
 We use the following notation for the entries to 
 matrices  $B_N^{-1}({d,q})$:
   \begin{multline}\label{Binv}
   B_N^{-1}({d,q})=\\=
   \left[
   \begin{array}{cccccc}
      \alpha_{N,0,N}(d,q)
     &\dots&
       \alpha_{N,0,1}(d,q)&
        \beta_{N,0,1}(d,q)
     &\dots&
       \beta_{N,0,N}(d,q)
     \\
     \vdots&\ddots&\vdots& \vdots&\ddots&\vdots\\
         \alpha_{N,2N-1,N}(d,q)
      &\dots&
       \alpha_{N,2N-1,1}(d,q)& \beta_{N,2N-1,1}(d,q)
      &\dots&
       \beta_{N,2N-1,N}(d,q)
   \end{array}
   \right].
 \end{multline}

In this paper  we analyse the top and the bottom rows 
of these inverse matrices only, that is, the  numbers
$ \alpha_{N,0,1}(d,q)$, \dots, $\alpha_{N,0,N}(d,q)$,
$ \beta_{N,0,1}(d,q)$,\dots, $\beta_{N,0,N}(d,q)$
and
$ \alpha_{N,2N-1,1}(d,q)$,\dots, $\alpha_{N,2N-1,N}(d,q)$,
$ \beta_{N,2N-1,1}(d,q)$,\dots, $\beta_{N,2N-1,N}(d,q)$.

 \section{Numerical observations and conjectures}
 \numberwithin{equation}{subsection}
 
 In this section we present some numerical 
 data and conjectures based on them; 
 other related numerical data can be found 
 at \url{https://doi.org/10.13140/RG.2.2.12349.22247}.

\subsection{Bottom row}

Numerical data suggest that the bottom
row of matrix $B_N^{-1}({d,q})$
is antisymmetrical.

\begin{myconjecture}\footnote{\emph{Added in version 2.} Andriy Bondarenko 
and Andrius Grigutis independently informed the author that  ChatGPT was able to prove this Conjecture.}
  For  $d=0$ or $d=1$, for all $q$,  $N$, and $j=1, \dots,N$,
  \begin{equation}
    \alpha_{N,2N-1,j}(d,q)= -\beta_{N,2N-1,j}(d,q).
  \end{equation}
\end{myconjecture}

Clearly, this conjecture can be restated as the equality of the
corresponding minors of the matrix $B_N({d,q})$.
For small values of $N$, the conjecture can be easily verified 
by any  computer algebra system.

 \subsection{Limiting values of 
 $\alpha_{N,0,j}(0,q)$ and $\beta_{N,0,j}(0,q)$}

In the rest of the paper, we consider
the entries to the top row of matrix
 $B_N^{-1}({0,q})$. In contrast to the bottom row, 
 there is no direct relationship between 
 the $\alpha$-entries (from the left-hand half
 of the row)
 and the $\beta$-entries (from the  right-hand half),
 thus 
 these numbers should be considered separately.

 Numerical calculations suggest the existence of 
 limiting values when $N$ tends to infinity over 
 an arithmetic progression modulo~$q$.

 \begin{myconjecture}\label{limexists}
   For every $q$, $M$, and $j$, there are (finite or infinite) limits
\begin{equation}\label{alims}
  {\alpha_{M,j}^\infty(q)=}
   \lim_{L\rightarrow\infty}\alpha_{qL+M,0,j}(0,q)  
   \qquad \textrm{and} \qquad  
{\beta_{M,j}^\infty(q)=}
  \lim_{L\rightarrow\infty}\beta_{qL+M,0,j}(0,q).
\end{equation}
\end{myconjecture}

It follows from the definitions \eqref{alims} 
that numbers $\alpha^\infty_{M,j}$
and  $\beta^\infty_{M,j}$
are periodic in $M$: 
  \begin{equation}\label{perM}
    \alpha_{M+q,j}^\infty(q)=\alpha_{M,j}^\infty(q)
    \qquad \textrm{and} \qquad \beta_{M+q,j}^\infty(q)=\beta_{M,j}^\infty(q).
  \end{equation}
  Numerical data indicate similar periodicity in $j$.
 \begin{myconjecture}\label{perij}
 For every $q$, $M$,   $j$, 
  \begin{equation}\label{perj}
    \alpha_{M,j+q}^\infty(q)=\alpha_{M,j}^\infty(q)
     \qquad and \qquad 
       \beta_{M,j+q}^\infty(q)=\beta_{M,j}^\infty(q).
  \end{equation}
  \end{myconjecture}

 Besides \eqref{perj}, the numerical data indicate 
 other equalities between numbers  
$\alpha_{M,j}^\infty(q)$ and $\beta_{M,j}^\infty(q)$.
  \begin{myconjecture}\label{consymm}
 For every $q$, $M$,   $j_1$, and $j_2$,
 if $j_1+j_2\equiv 0\bmod{q}$ then
  \begin{equation}\label{symj}
    \alpha_{M,j_1}^\infty(q)=\alpha_{M,j_2}^\infty(q)
     \qquad and \qquad 
       \beta_{M,j_1}^\infty(q)=\beta_{M,j_2}^\infty(q).
  \end{equation}
  \end{myconjecture}

   The numerical data allow us to state a number of conjectures about the concrete values of   $\alpha_{M,0}^\infty(q)$ and $\beta_{M,0}^\infty(q)$.
 
   \begin{table}[h] \center
 
\newlength{\cnqsnukiq}
\setlength{\cnqsnukiq}{\widthof{}}
\setlength{\cnqsnukiq}{\maxof{%
  \cnqsnukiq}{\widthof{$2$}}}
\setlength{\cnqsnukiq}{\maxof{%
  \cnqsnukiq}{\widthof{$3$}}}
\setlength{\cnqsnukiq}{\maxof{%
  \cnqsnukiq}{\widthof{$4$}}}
\setlength{\cnqsnukiq}{\maxof{%
  \cnqsnukiq}{\widthof{$5$}}}
\setlength{\cnqsnukiq}{\maxof{%
  \cnqsnukiq}{\widthof{$6$}}}
\setlength{\cnqsnukiq}{\maxof{%
  \cnqsnukiq}{\widthof{$7$}}}
\setlength{\cnqsnukiq}{\maxof{%
  \cnqsnukiq}{\widthof{$8$}}}
\setlength{\cnqsnukiq}{\maxof{%
  \cnqsnukiq}{\widthof{$9$}}}
\setlength{\cnqsnukiq}{\maxof{%
  \cnqsnukiq}{\widthof{$10$}}}
\setlength{\cnqsnukiq}{\maxof{%
  \cnqsnukiq}{\widthof{$11$}}}
\setlength{\cnqsnukiq}{\maxof{%
  \cnqsnukiq}{\widthof{$12$}}}
\newlength{\tbopnsfdq}
\setlength{\tbopnsfdq}{\widthof{}}
\setlength{\tbopnsfdq}{\maxof{%
  \tbopnsfdq}{\widthof{$-2-6.46\dotsc\cdot 10^{-271}$}}}
\setlength{\tbopnsfdq}{\maxof{%
  \tbopnsfdq}{\widthof{$-2+4.95\dotsc\cdot 10^{-169}$}}}
\setlength{\tbopnsfdq}{\maxof{%
  \tbopnsfdq}{\widthof{$-2+5.89\dotsc\cdot 10^{-120}$}}}
\setlength{\tbopnsfdq}{\maxof{%
  \tbopnsfdq}{\widthof{$-2-4.40\dotsc\cdot 10^{-89}$}}}
\setlength{\tbopnsfdq}{\maxof{%
  \tbopnsfdq}{\widthof{$-2+3.99\dotsc\cdot 10^{-79}$}}}
\setlength{\tbopnsfdq}{\maxof{%
  \tbopnsfdq}{\widthof{$-2+1.73\dotsc\cdot 10^{-55}$}}}
\setlength{\tbopnsfdq}{\maxof{%
  \tbopnsfdq}{\widthof{$-2-3.69\dotsc\cdot 10^{-50}$}}}
\setlength{\tbopnsfdq}{\maxof{%
  \tbopnsfdq}{\widthof{$-2-7.08\dotsc\cdot 10^{-37}$}}}
\setlength{\tbopnsfdq}{\maxof{%
  \tbopnsfdq}{\widthof{$-2+2.69\dotsc\cdot 10^{-39}$}}}
\setlength{\tbopnsfdq}{\maxof{%
  \tbopnsfdq}{\widthof{$-2-4.37\dotsc\cdot 10^{-28}$}}}
\setlength{\tbopnsfdq}{\maxof{%
  \tbopnsfdq}{\widthof{$-2+5.85\dotsc\cdot 10^{-27}$}}}
\newlength{\pedbydmon}
\setlength{\pedbydmon}{\widthof{}}
\setlength{\pedbydmon}{\maxof{%
  \pedbydmon}{\widthof{$-2-1.44\dotsc\cdot 10^{-1085}$}}}
\setlength{\pedbydmon}{\maxof{%
  \pedbydmon}{\widthof{$-2+7.31\dotsc\cdot 10^{-702}$}}}
\setlength{\pedbydmon}{\maxof{%
  \pedbydmon}{\widthof{$-2+1.12\dotsc\cdot 10^{-511}$}}}
\setlength{\pedbydmon}{\maxof{%
  \pedbydmon}{\widthof{$-2-6.71\dotsc\cdot 10^{-396}$}}}
\setlength{\pedbydmon}{\maxof{%
  \pedbydmon}{\widthof{$-2+7.50\dotsc\cdot 10^{-338}$}}}
\setlength{\pedbydmon}{\maxof{%
  \pedbydmon}{\widthof{$-2-1.56\dotsc\cdot 10^{-279}$}}}
\setlength{\pedbydmon}{\maxof{%
  \pedbydmon}{\widthof{$-2-3.06\dotsc\cdot 10^{-237}$}}}
\setlength{\pedbydmon}{\maxof{%
  \pedbydmon}{\widthof{$-2+1.28\dotsc\cdot 10^{-203}$}}}
\setlength{\pedbydmon}{\maxof{%
  \pedbydmon}{\widthof{$-2-7.83\dotsc\cdot 10^{-188}$}}}
\setlength{\pedbydmon}{\maxof{%
  \pedbydmon}{\widthof{$-2+1.10\dotsc\cdot 10^{-162}$}}}
\setlength{\pedbydmon}{\maxof{%
  \pedbydmon}{\widthof{$-2+1.03\dotsc\cdot 10^{-153}$}}}
\newlength{\jmsokjqpb}
\setlength{\jmsokjqpb}{\widthof{}}
\setlength{\jmsokjqpb}{\maxof{%
  \jmsokjqpb}{\widthof{$-2-1.10\dotsc\cdot 10^{-4351}$}}}
\setlength{\jmsokjqpb}{\maxof{%
  \jmsokjqpb}{\widthof{$-2+1.50\dotsc\cdot 10^{-2857}$}}}
\setlength{\jmsokjqpb}{\maxof{%
  \jmsokjqpb}{\widthof{$-2+2.88\dotsc\cdot 10^{-2113}$}}}
\setlength{\jmsokjqpb}{\maxof{%
  \jmsokjqpb}{\widthof{$-2-2.83\dotsc\cdot 10^{-1664}$}}}
\setlength{\jmsokjqpb}{\maxof{%
  \jmsokjqpb}{\widthof{$-2+9.55\dotsc\cdot 10^{-1404}$}}}
\setlength{\jmsokjqpb}{\maxof{%
  \jmsokjqpb}{\widthof{$-2-7.61\dotsc\cdot 10^{-1184}$}}}
\setlength{\jmsokjqpb}{\maxof{%
  \jmsokjqpb}{\widthof{$-2-3.26\dotsc\cdot 10^{-1020}$}}}
\setlength{\jmsokjqpb}{\maxof{%
  \jmsokjqpb}{\widthof{$-2+3.70\dotsc\cdot 10^{-870}$}}}
\setlength{\jmsokjqpb}{\maxof{%
  \jmsokjqpb}{\widthof{$-2-2.03\dotsc\cdot 10^{-812}$}}}
\setlength{\jmsokjqpb}{\maxof{%
  \jmsokjqpb}{\widthof{$-2+3.18\dotsc\cdot 10^{-763}$}}}
\setlength{\jmsokjqpb}{\maxof{%
  \jmsokjqpb}{\widthof{$-2+6.47\dotsc\cdot 10^{-656}$}}}
\begin{tabular}{||r||c|c|c||}
\hhline{|t:=:t:===:t|}
  \multicolumn{1}{||r||}{$q$}
&
  \multicolumn{1}{c|}{$N=20$}
&
  \multicolumn{1}{c|}{$N=40$}
&
  \multicolumn{1}{c||}{$N=80$}
\\
\hhline{||-||-|-|-||}

\parbox{\cnqsnukiq}{$ 2 $}

&\parbox{\tbopnsfdq}{$-2-6.46\dotsc\cdot 10^{-271}$}

&\parbox{\pedbydmon}{$-2-1.44\dotsc\cdot 10^{-1085}$}

&\parbox{\jmsokjqpb}{$-2-1.10\dotsc\cdot 10^{-4351}$}
\\ 
\parbox{\cnqsnukiq}{$ 3 $}

&\parbox{\tbopnsfdq}{$-2+4.95\dotsc\cdot 10^{-169}$}

&\parbox{\pedbydmon}{$-2+7.31\dotsc\cdot 10^{-702}$}

&\parbox{\jmsokjqpb}{$-2+1.50\dotsc\cdot 10^{-2857}$}
\\ 
\parbox{\cnqsnukiq}{$ 4 $}

&\parbox{\tbopnsfdq}{$-2+5.89\dotsc\cdot 10^{-120}$}

&\parbox{\pedbydmon}{$-2+1.12\dotsc\cdot 10^{-511}$}

&\parbox{\jmsokjqpb}{$-2+2.88\dotsc\cdot 10^{-2113}$}
\\ 
\parbox{\cnqsnukiq}{$ 5 $}

&\parbox{\tbopnsfdq}{$-2-4.40\dotsc\cdot 10^{-89}$}

&\parbox{\pedbydmon}{$-2-6.71\dotsc\cdot 10^{-396}$}

&\parbox{\jmsokjqpb}{$-2-2.83\dotsc\cdot 10^{-1664}$}
\\ 
\parbox{\cnqsnukiq}{$ 6 $}

&\parbox{\tbopnsfdq}{$-2+3.99\dotsc\cdot 10^{-79}$}

&\parbox{\pedbydmon}{$-2+7.50\dotsc\cdot 10^{-338}$}

&\parbox{\jmsokjqpb}{$-2+9.55\dotsc\cdot 10^{-1404}$}
\\ 
\parbox{\cnqsnukiq}{$ 7 $}

&\parbox{\tbopnsfdq}{$-2+1.73\dotsc\cdot 10^{-55}$}

&\parbox{\pedbydmon}{$-2-1.56\dotsc\cdot 10^{-279}$}

&\parbox{\jmsokjqpb}{$-2-7.61\dotsc\cdot 10^{-1184}$}
\\ 
\parbox{\cnqsnukiq}{$ 8 $}

&\parbox{\tbopnsfdq}{$-2-3.69\dotsc\cdot 10^{-50}$}

&\parbox{\pedbydmon}{$-2-3.06\dotsc\cdot 10^{-237}$}

&\parbox{\jmsokjqpb}{$-2-3.26\dotsc\cdot 10^{-1020}$}
\\ 
\parbox{\cnqsnukiq}{$ 9 $}

&\parbox{\tbopnsfdq}{$-2-7.08\dotsc\cdot 10^{-37}$}

&\parbox{\pedbydmon}{$-2+1.28\dotsc\cdot 10^{-203}$}

&\parbox{\jmsokjqpb}{$-2+3.70\dotsc\cdot 10^{-870}$}
\\ 
\parbox{\cnqsnukiq}{$ 10 $}

&\parbox{\tbopnsfdq}{$-2+2.69\dotsc\cdot 10^{-39}$}

&\parbox{\pedbydmon}{$-2-7.83\dotsc\cdot 10^{-188}$}

&\parbox{\jmsokjqpb}{$-2-2.03\dotsc\cdot 10^{-812}$}
\\ 
\parbox{\cnqsnukiq}{$ 11 $}

&\parbox{\tbopnsfdq}{$-2-4.37\dotsc\cdot 10^{-28}$}

&\parbox{\pedbydmon}{$-2+1.10\dotsc\cdot 10^{-162}$}

&\parbox{\jmsokjqpb}{$-2+3.18\dotsc\cdot 10^{-763}$}
\\ 
\parbox{\cnqsnukiq}{$ 12 $}

&\parbox{\tbopnsfdq}{$-2+5.85\dotsc\cdot 10^{-27}$}

&\parbox{\pedbydmon}{$-2+1.03\dotsc\cdot 10^{-153}$}

&\parbox{\jmsokjqpb}{$-2+6.47\dotsc\cdot 10^{-656}$}
\\ 
\hhline{|b:=:b:===:b|}
\end{tabular}
\nocaption{Values of $\alpha_{N,0,q}(0,q)$}
 
\caption{Values of $\alpha_{N,0,q}(0,q)$}
\label{tabajis0}
\end{table}
 
    \begin{table}[h] \center
 
\newlength{\ioaviwcmo}
\setlength{\ioaviwcmo}{\widthof{}}
\setlength{\ioaviwcmo}{\maxof{%
  \ioaviwcmo}{\widthof{$2$}}}
\setlength{\ioaviwcmo}{\maxof{%
  \ioaviwcmo}{\widthof{$3$}}}
\setlength{\ioaviwcmo}{\maxof{%
  \ioaviwcmo}{\widthof{$4$}}}
\setlength{\ioaviwcmo}{\maxof{%
  \ioaviwcmo}{\widthof{$5$}}}
\setlength{\ioaviwcmo}{\maxof{%
  \ioaviwcmo}{\widthof{$6$}}}
\setlength{\ioaviwcmo}{\maxof{%
  \ioaviwcmo}{\widthof{$7$}}}
\setlength{\ioaviwcmo}{\maxof{%
  \ioaviwcmo}{\widthof{$8$}}}
\setlength{\ioaviwcmo}{\maxof{%
  \ioaviwcmo}{\widthof{$9$}}}
\setlength{\ioaviwcmo}{\maxof{%
  \ioaviwcmo}{\widthof{$10$}}}
\setlength{\ioaviwcmo}{\maxof{%
  \ioaviwcmo}{\widthof{$11$}}}
\setlength{\ioaviwcmo}{\maxof{%
  \ioaviwcmo}{\widthof{$12$}}}
\newlength{\kolwxhxbg}
\setlength{\kolwxhxbg}{\widthof{}}
\setlength{\kolwxhxbg}{\maxof{%
  \kolwxhxbg}{\widthof{$\phantom{-}6.74\dotsc\cdot 10^{-271}$}}}
\setlength{\kolwxhxbg}{\maxof{%
  \kolwxhxbg}{\widthof{$-4.88\dotsc\cdot 10^{-169}$}}}
\setlength{\kolwxhxbg}{\maxof{%
  \kolwxhxbg}{\widthof{$\phantom{-}5.89\dotsc\cdot 10^{-120}$}}}
\setlength{\kolwxhxbg}{\maxof{%
  \kolwxhxbg}{\widthof{$-4.40\dotsc\cdot 10^{-89}$}}}
\setlength{\kolwxhxbg}{\maxof{%
  \kolwxhxbg}{\widthof{$-7.36\dotsc\cdot 10^{-79}$}}}
\setlength{\kolwxhxbg}{\maxof{%
  \kolwxhxbg}{\widthof{$-1.59\dotsc\cdot 10^{-55}$}}}
\setlength{\kolwxhxbg}{\maxof{%
  \kolwxhxbg}{\widthof{$-3.69\dotsc\cdot 10^{-50}$}}}
\setlength{\kolwxhxbg}{\maxof{%
  \kolwxhxbg}{\widthof{$-7.08\dotsc\cdot 10^{-37}$}}}
\setlength{\kolwxhxbg}{\maxof{%
  \kolwxhxbg}{\widthof{$\phantom{-}1.90\dotsc\cdot 10^{-38}$}}}
\setlength{\kolwxhxbg}{\maxof{%
  \kolwxhxbg}{\widthof{$\phantom{-}3.42\dotsc\cdot 10^{-28}$}}}
\setlength{\kolwxhxbg}{\maxof{%
  \kolwxhxbg}{\widthof{$\phantom{-}5.85\dotsc\cdot 10^{-27}$}}}
\newlength{\fxwbvzzbj}
\setlength{\fxwbvzzbj}{\widthof{}}
\setlength{\fxwbvzzbj}{\maxof{%
  \fxwbvzzbj}{\widthof{$\phantom{-}1.45\dotsc\cdot 10^{-1085}$}}}
\setlength{\fxwbvzzbj}{\maxof{%
  \fxwbvzzbj}{\widthof{$-7.40\dotsc\cdot 10^{-702}$}}}
\setlength{\fxwbvzzbj}{\maxof{%
  \fxwbvzzbj}{\widthof{$\phantom{-}1.12\dotsc\cdot 10^{-511}$}}}
\setlength{\fxwbvzzbj}{\maxof{%
  \fxwbvzzbj}{\widthof{$-6.71\dotsc\cdot 10^{-396}$}}}
\setlength{\fxwbvzzbj}{\maxof{%
  \fxwbvzzbj}{\widthof{$-7.39\dotsc\cdot 10^{-338}$}}}
\setlength{\fxwbvzzbj}{\maxof{%
  \fxwbvzzbj}{\widthof{$\phantom{-}1.54\dotsc\cdot 10^{-279}$}}}
\setlength{\fxwbvzzbj}{\maxof{%
  \fxwbvzzbj}{\widthof{$-3.06\dotsc\cdot 10^{-237}$}}}
\setlength{\fxwbvzzbj}{\maxof{%
  \fxwbvzzbj}{\widthof{$\phantom{-}1.28\dotsc\cdot 10^{-203}$}}}
\setlength{\fxwbvzzbj}{\maxof{%
  \fxwbvzzbj}{\widthof{$\phantom{-}1.42\dotsc\cdot 10^{-187}$}}}
\setlength{\fxwbvzzbj}{\maxof{%
  \fxwbvzzbj}{\widthof{$-1.11\dotsc\cdot 10^{-162}$}}}
\setlength{\fxwbvzzbj}{\maxof{%
  \fxwbvzzbj}{\widthof{$\phantom{-}1.03\dotsc\cdot 10^{-153}$}}}
\newlength{\lgvawusjj}
\setlength{\lgvawusjj}{\widthof{}}
\setlength{\lgvawusjj}{\maxof{%
  \lgvawusjj}{\widthof{$\phantom{-}1.10\dotsc\cdot 10^{-4351}$}}}
\setlength{\lgvawusjj}{\maxof{%
  \lgvawusjj}{\widthof{$-1.50\dotsc\cdot 10^{-2857}$}}}
\setlength{\lgvawusjj}{\maxof{%
  \lgvawusjj}{\widthof{$\phantom{-}2.88\dotsc\cdot 10^{-2113}$}}}
\setlength{\lgvawusjj}{\maxof{%
  \lgvawusjj}{\widthof{$-2.83\dotsc\cdot 10^{-1664}$}}}
\setlength{\lgvawusjj}{\maxof{%
  \lgvawusjj}{\widthof{$-9.88\dotsc\cdot 10^{-1404}$}}}
\setlength{\lgvawusjj}{\maxof{%
  \lgvawusjj}{\widthof{$\phantom{-}7.90\dotsc\cdot 10^{-1184}$}}}
\setlength{\lgvawusjj}{\maxof{%
  \lgvawusjj}{\widthof{$-3.26\dotsc\cdot 10^{-1020}$}}}
\setlength{\lgvawusjj}{\maxof{%
  \lgvawusjj}{\widthof{$\phantom{-}3.69\dotsc\cdot 10^{-870}$}}}
\setlength{\lgvawusjj}{\maxof{%
  \lgvawusjj}{\widthof{$\phantom{-}2.32\dotsc\cdot 10^{-812}$}}}
\setlength{\lgvawusjj}{\maxof{%
  \lgvawusjj}{\widthof{$-3.66\dotsc\cdot 10^{-763}$}}}
\setlength{\lgvawusjj}{\maxof{%
  \lgvawusjj}{\widthof{$\phantom{-}6.47\dotsc\cdot 10^{-656}$}}}
\begin{tabular}{||r||c|c|c||}
\hhline{|t:=:t:===:t|}
  \multicolumn{1}{||r||}{$q$}
&
  \multicolumn{1}{c|}{$N=20$}
&
  \multicolumn{1}{c|}{$N=40$}
&
  \multicolumn{1}{c||}{$N=80$}
\\
\hhline{||-||-|-|-||}

\parbox{\ioaviwcmo}{$ 2 $}

&\parbox{\kolwxhxbg}{$ \phantom{-}6.74\dotsc\cdot 10^{-271}$}

&\parbox{\fxwbvzzbj}{$ \phantom{-}1.45\dotsc\cdot 10^{-1085}$}

&\parbox{\lgvawusjj}{$ \phantom{-}1.10\dotsc\cdot 10^{-4351}$}
\\ 
\parbox{\ioaviwcmo}{$ 3 $}

&\parbox{\kolwxhxbg}{$ -4.88\dotsc\cdot 10^{-169}$}

&\parbox{\fxwbvzzbj}{$ -7.40\dotsc\cdot 10^{-702}$}

&\parbox{\lgvawusjj}{$ -1.50\dotsc\cdot 10^{-2857}$}
\\ 
\parbox{\ioaviwcmo}{$ 4 $}

&\parbox{\kolwxhxbg}{$ \phantom{-}5.89\dotsc\cdot 10^{-120}$}

&\parbox{\fxwbvzzbj}{$ \phantom{-}1.12\dotsc\cdot 10^{-511}$}

&\parbox{\lgvawusjj}{$ \phantom{-}2.88\dotsc\cdot 10^{-2113}$}
\\ 
\parbox{\ioaviwcmo}{$ 5 $}

&\parbox{\kolwxhxbg}{$ -4.40\dotsc\cdot 10^{-89}$}

&\parbox{\fxwbvzzbj}{$ -6.71\dotsc\cdot 10^{-396}$}

&\parbox{\lgvawusjj}{$ -2.83\dotsc\cdot 10^{-1664}$}
\\ 
\parbox{\ioaviwcmo}{$ 6 $}

&\parbox{\kolwxhxbg}{$ -7.36\dotsc\cdot 10^{-79}$}

&\parbox{\fxwbvzzbj}{$ -7.39\dotsc\cdot 10^{-338}$}

&\parbox{\lgvawusjj}{$ -9.88\dotsc\cdot 10^{-1404}$}
\\ 
\parbox{\ioaviwcmo}{$ 7 $}

&\parbox{\kolwxhxbg}{$ -1.59\dotsc\cdot 10^{-55}$}

&\parbox{\fxwbvzzbj}{$ \phantom{-}1.54\dotsc\cdot 10^{-279}$}

&\parbox{\lgvawusjj}{$ \phantom{-}7.90\dotsc\cdot 10^{-1184}$}
\\ 
\parbox{\ioaviwcmo}{$ 8 $}

&\parbox{\kolwxhxbg}{$ -3.69\dotsc\cdot 10^{-50}$}

&\parbox{\fxwbvzzbj}{$ -3.06\dotsc\cdot 10^{-237}$}

&\parbox{\lgvawusjj}{$ -3.26\dotsc\cdot 10^{-1020}$}
\\ 
\parbox{\ioaviwcmo}{$ 9 $}

&\parbox{\kolwxhxbg}{$ -7.08\dotsc\cdot 10^{-37}$}

&\parbox{\fxwbvzzbj}{$ \phantom{-}1.28\dotsc\cdot 10^{-203}$}

&\parbox{\lgvawusjj}{$ \phantom{-}3.69\dotsc\cdot 10^{-870}$}
\\ 
\parbox{\ioaviwcmo}{$ 10 $}

&\parbox{\kolwxhxbg}{$ \phantom{-}1.90\dotsc\cdot 10^{-38}$}

&\parbox{\fxwbvzzbj}{$ \phantom{-}1.42\dotsc\cdot 10^{-187}$}

&\parbox{\lgvawusjj}{$ \phantom{-}2.32\dotsc\cdot 10^{-812}$}
\\ 
\parbox{\ioaviwcmo}{$ 11 $}

&\parbox{\kolwxhxbg}{$ \phantom{-}3.42\dotsc\cdot 10^{-28}$}

&\parbox{\fxwbvzzbj}{$ -1.11\dotsc\cdot 10^{-162}$}

&\parbox{\lgvawusjj}{$ -3.66\dotsc\cdot 10^{-763}$}
\\ 
\parbox{\ioaviwcmo}{$ 12 $}

&\parbox{\kolwxhxbg}{$ \phantom{-}5.85\dotsc\cdot 10^{-27}$}

&\parbox{\fxwbvzzbj}{$ \phantom{-}1.03\dotsc\cdot 10^{-153}$}

&\parbox{\lgvawusjj}{$ \phantom{-}6.47\dotsc\cdot 10^{-656}$}
\\ 
\hhline{|b:=:b:===:b|}
\end{tabular}
\nocaption{Values of $\beta_{N,0,q}(0,q)$}
 
\caption{Values of $\beta_{N,0,q}(0,q)$}
\label{tabbjis0}
\end{table}
 
\begin{myconjecture}\label{j=0}
  For all $q$,  $M$,  and $j$, if $j\equiv 0\bmod{q}$
  then
\begin{equation}\label{20}
\alpha^\infty_{M,j}(q)=-2
      \qquad \textrm{and} \qquad 
      \beta^\infty_{M,j}(q)=0.
   \end{equation} 
\end{myconjecture}    
\noindent
Tables \ref{tabajis0}--\ref{tabbjis0} and illustrates this conjecture.
  
  As for other finite  values of  $\alpha^\infty_{M,j}(q)$
  and $\beta^\infty_{M,j}(q)$,
  it seems that all of them are algebraic number.
  \begin{myconjecture}
  For all $q$, $j$, and $M$,
  if the  value of  $\alpha^\infty_{M,j}(q)$ 
  is finite, then it
 belongs to the totally 
  real cyclotomic 
  field~$\mathbb{Q}(\cos(\pi/q))$; moreover, if $q$ is even,
  then  $\alpha^\infty_{M,j}(q)$ belongs already to~$\mathbb{Q}(\cos(2\pi/q))$.
\end{myconjecture}    

 Finite $\beta^\infty_{M,j}(q)$ belongs to the same cyclotomic
 field, but after the scaling by $\sqrt{q}$. 

\begin{myconjecture}
  For all $q$, $j$, and $M$,
    if the  value  of  $\beta^\infty_{M,j}(q)$ is finite, then 
 $\sqrt{q}\beta^\infty_{M,j}(q)$
 belongs to the totally 
  real cyclotomic 
  field~$\mathbb{Q}(\cos(\pi/q))$; moreover, if $q$ is even,
  then  $\sqrt{q}\beta^\infty_{M,j}(q)$ belongs already to its subfield
  ~$\mathbb{Q}(\cos(2\pi/q))$.
\end{myconjecture}

  In Tables \ref{tab2alphaInfMj}--\ref{tab12betaInfMj}, we
  present our conjectures about the  values of  $\alpha^\infty_{M,j}(q)$ and $\sqrt{q}\beta^\infty_{M,j}(q)$
for $q=2,\dots,12$, $ M=0, \dots, \  q-1$ and 
$j=1,\ \dots.\  q-1$. 
 In these tables $c_q=\cos(\pi/q)$. The tables allow one to
 determine the expected values of  $\alpha^\infty_{M,j}(q)$ and $\sqrt{q}\beta^\infty_{M,j}(q)$  
 for all  other values of  $M$ and $j$ via
    \eqref{perM}, \eqref{perj}, and \eqref{20}.
   Tables \ref{tabalphalim} and \ref{tabbetalim} demonstrate the rate of convergence 
   of   $\alpha_{M,0,j}(0,q)$ and $\sqrt{q}\beta_{M,0,j}(0,q)$
   to the expected finite limiting values.

 The appearance of the field ~$\mathbb{Q}(\cos(\pi/q))$ and above-described periodicities
 modulo~$q$ are rather surprising by the following reason:
 number $q$ 
 occurs in the definition of matrices  $B_{N}(0,q)$
 in the denominators of certain real-valued expressions
 only. Hypothetically, the periodicities could be connected 
  to the periodicity of the characters, but the latter  are not used in the definition \eqref{defB} of matrix $B_{N}(0,q)$ at all.
 
 We were able to state a number of conjecture about  
 numbers  $\alpha^\infty_{M,j}(q)$ and $\beta^\infty_{M,j}(q)$;
 however, no compact unifying description
  was found for Tables~\ref{tab2alphaInfMj}--\ref{tab12betaInfMj}
  and similar tables for larger~$q$.

\begin{table}[h] \center

\begin{tabular}{||c||c||}
\hhline{|t:=:t:=:t|}
  \multicolumn{1}{||c||}{$M$}
&
  \multicolumn{1}{c||}{$j=1$}
\\
\hhline{||-||-||}

0, 1
&{$ 2 $}
\\ 
\hhline{|b:=:b:=:b|}
\end{tabular}
\nocaption{Expected values of $\alpha^\infty_{M,j}(2)$}
 
\caption{Expected values $\alpha^\infty_{M,j}(2)$}
\label{tab2alphaInfMj}
\end{table}  
 
\begin{table}[h] \center

\begin{tabular}{||c||c||}
\hhline{|t:=:t:=:t|}
  \multicolumn{1}{||c||}{$M$}
&
  \multicolumn{1}{c||}{$j=1$}
\\
\hhline{||-||-||}

0, 1
&{$ 4 $}
\\ 
\hhline{|b:=:b:=:b|}
\end{tabular}
\nocaption{Expected values of $\sqrt{2}\beta^\infty_{M,j}(2)$}
 
\caption{Expected values of $\sqrt{2}\beta^\infty_{M,j}(2)$}
\label{tab2betaInfMj}
\end{table}

\begin{table}[h] \center

\begin{tabular}{||c||c||}
\hhline{|t:=:t:=:t|}
  \multicolumn{1}{||c||}{$M$}
&
  \multicolumn{1}{c||}{$j=1,\,2$}
\\
\hhline{||-||-||}

0, 1, 2
&$1$\\ 
\hhline{|b:=:b:=:b|}
\end{tabular}
 
\caption{Expected values of $\alpha^\infty_{M,j}(3)$}
\label{tab3alphaInfMj}
\end{table}  
 
\begin{table}[h] \center

\begin{tabular}{||c||c||}
\hhline{|t:=:t:=:t|}
  \multicolumn{1}{||c||}{$M$}
&
  \multicolumn{1}{c||}{$j=1,\,2$}
\\
\hhline{||-||-||}

0, 1, 2
&$3$\\ 
\hhline{|b:=:b:=:b|}
\end{tabular}
 
\caption{Expected values of $\sqrt{3}\beta^\infty_{M,j}(3)$}
\label{tab3betaInfMj}
\end{table}

\begin{table}[h] \center

\begin{tabular}{||c||c|c||}
\hhline{|t:=:t:==:t|}
  \multicolumn{1}{||c||}{$M$}
&
  \multicolumn{1}{c|}{$j=1,\,3$}
&
  \multicolumn{1}{c||}{$j=2$}
\\
\hhline{||-||-|-||}

0, 2
&{$ 1 $}

&{$ 0 $}
\\ 
1
&{$ 1/2 $}

&{$ 1 $}
\\ 
\hhline{|b:=:b:==:b|}
\end{tabular}
\nocaption{Expected values of $\alpha^\infty_{M,j}(4)$}
 
\caption{Expected values of $\alpha^\infty_{M,j}(4)$}
\label{tab4alphaInfMj}
\end{table}  

\begin{table}[h] \center

\begin{tabular}{||c||c|c||}
\hhline{|t:=:t:==:t|}
  \multicolumn{1}{||c||}{$M$}
&
  \multicolumn{1}{c|}{$j=1,\,3$}
&
  \multicolumn{1}{c||}{$j=2$}
\\
\hhline{||-||-|-||}

0, 2
&{$ 2 $}

&{$ 4 $}
\\ 
1
&{$ 3 $}

&{$ 2 $}
\\ 
\hhline{|b:=:b:==:b|}
\end{tabular}
\nocaption{Expected values of $\sqrt{4}\beta^\infty_{M,j}(4)$}
 
\caption{Expected values of $\sqrt{4}\beta^\infty_{M,j}(4)$}
\label{tab4betaInfMj}
\end{table}

\begin{table}[h] \center

\begin{tabular}{||c||c|c||}
\hhline{|t:=:t:==:t|}
  \multicolumn{1}{||c||}{$M$}
&
  \multicolumn{1}{c|}{$j=1,\,4$}
&
  \multicolumn{1}{c||}{$j=2,\,3$}
\\
\hhline{||-||-|-||}

0, 3
&$c_5$
&$-c_5+1$\\ 
1, 2
&$-c_5+1$
&$c_5$\\ 
4
&$+\infty$
&$-\infty$\\ 
\hhline{|b:=:b:==:b|}
\end{tabular}
 
\caption{Expected values of $\alpha^\infty_{M,j}(5)$}
\label{tab5alphaInfMj}
\end{table}  

\begin{table}[h] \center

\begin{tabular}{||c||c|c||}
\hhline{|t:=:t:==:t|}
  \multicolumn{1}{||c||}{$M$}
&
  \multicolumn{1}{c|}{$j=1,\,4$}
&
  \multicolumn{1}{c||}{$j=2,\,3$}
\\
\hhline{||-||-|-||}

0, 3
&$c_5+1$
&$-c_5+4$\\ 
1, 2
&$-c_5+4$
&$c_5+1$\\ 
4
&$-\infty$
&$+\infty$\\ 
\hhline{|b:=:b:==:b|}
\end{tabular}
 
\caption{Expected values of $\sqrt{5}\beta^\infty_{M,j}(5)$}
\label{tab5betaInfMj}
\end{table}   

\begin{table}[h] \center

\begin{tabular}{||c||c|c|c||}
\hhline{|t:=:t:===:t|}
  \multicolumn{1}{||c||}{$M$}
&
  \multicolumn{1}{c|}{$j=1,\,5$}
&
  \multicolumn{1}{c|}{$j=2,\,4$}
&
  \multicolumn{1}{c||}{$j=3$}
\\
\hhline{||-||-|-|-||}

0, 4
&$0$
&$4$
&$-6$\\ 
1, 3
&$0$
&$0$
&$2$\\ 
2
&$0$
&$1$
&$0$\\ 
5
&$+\infty$
&$-\infty$
&$+\infty$\\ 
\hhline{|b:=:b:===:b|}
\end{tabular}
\nocaption{Expected values of $\alpha^\infty_{M,j}(6)$}
 
\caption{Expected values of $\alpha^\infty_{M,j}(6)$}
\label{tab6alphaInfMj}
\end{table}  

\begin{table}[h] \center

\begin{tabular}{||c||c|c|c||}
\hhline{|t:=:t:===:t|}
  \multicolumn{1}{||c||}{$M$}
&
  \multicolumn{1}{c|}{$j=1,\,5$}
&
  \multicolumn{1}{c|}{$j=2,\,4$}
&
  \multicolumn{1}{c||}{$j=3$}
\\
\hhline{||-||-|-|-||}

0, 4
&$0$
&$12$
&$-12$\\ 
1, 3
&$4$
&$0$
&$4$\\ 
2
&$3$
&$3$
&$0$\\ 
5
&$+\infty$
&$-\infty$
&$+\infty$\\ 
\hhline{|b:=:b:===:b|}
\end{tabular}
\nocaption{Expected values of $\sqrt{6}\beta^\infty_{M,j}(6)$}

\caption{Expected values of $\sqrt{6}\beta^\infty_{M,j}(6)$} 
\label{tab6betaInfMj}
\end{table}

\begin{table}[h] \center
 
\begin{tabular}{||c||c|c|c||}
\hhline{|t:=:t:===:t|}
  \multicolumn{1}{||c||}{$M$}
&
  \multicolumn{1}{c|}{$j=1,\,6$}
&
  \multicolumn{1}{c|}{$j=2,\,5$}
&
  \multicolumn{1}{c||}{$j=3,\,4$}
\\
\hhline{||-||-|-|-||}

0, 5
&$0$
&$2c_7+1$
&$-2c_7$\\ 
1, 4
&$-4c_7^2+2c_7+1$
&$0$
&$4c_7^2-2c_7$\\ 
2, 3
&$-4c_7^2+3$
&$4c_7^2-2$
&$0$\\ 
6
&$+\infty$
&$-\infty$
&$+\infty$\\ 
\hhline{|b:=:b:===:b|}
\end{tabular}
 
\caption{Expected values of $\alpha^\infty_{M,j}(7)$}
\label{tab7alphaInfMj}
\end{table}  

\begin{table}[h] \center

\begin{tabular}{||c||c|c|c||}
\hhline{|t:=:t:===:t|}
  \multicolumn{1}{||c||}{$M$}
&
  \multicolumn{1}{c|}{$j=1,\,6$}
&
  \multicolumn{1}{c|}{$j=2,\,5$}
&
  \multicolumn{1}{c||}{$j=3,\,4$}
\\
\hhline{||-||-|-|-||}

0, 5
&$0$
&$8c_7^2+2c_7+1$
&$-8c_7^2-2c_7+6$\\ 
1, 4
&$-4c_7^2+6c_7+3$
&$0$
&$4c_7^2-6c_7+4$\\ 
2, 3
&$-12c_7^2+4c_7+9$
&$12c_7^2-4c_7-2$
&$0$\\ 
6
&$+\infty$
&$-\infty$
&$+\infty$\\ 
\hhline{|b:=:b:===:b|}
\end{tabular}
 
\caption{Expected values of $\sqrt{7}\beta^\infty_{M,j}(7)$}
\label{tab7betaInfMj}
\end{table}

\begin{table}[h] \center
  
\begin{tabular}{||c||c|c|c|c||}
\hhline{|t:=:t:====:t|}
  \multicolumn{1}{||c||}{$M$}
&
  \multicolumn{1}{c|}{$j=1,\,7$}
&
  \multicolumn{1}{c|}{$j=2,\,6$}
&
  \multicolumn{1}{c|}{$j=3,\,5$}
&
  \multicolumn{1}{c||}{$j=4$}
\\
\hhline{||-||-|-|-|-||}

0
&$0$
&$4c_8^2-1$
&$-4c_8^2+2$
&$0$\\ 
1
&$-2c_8^2+1$
&$0$
&$2c_8^2-1$
&$2$\\ 
2
&$-2c_8^2+1$
&$2c_8^2+\frac{1}{2}$
&$0$
&$-1$\\ 
3
&$-c_8^2+\frac{1}{2}$
&$1$
&$c_8^2-\frac{1}{2}$
&$0$\\ 
4
&$-4c_8^2+2$
&$4c_8^2-1$
&$0$
&$0$\\ 
5
&$2c_8^2-1$
&$0$
&$-2c_8^2+1$
&$2$\\ 
6
&$0$
&$+\infty$
&$-\infty$
&$+\infty$\\ 
7
&$+\infty$
&$-\infty$
&$+\infty$
&$-\infty$\\ 
\hhline{|b:=:b:====:b|}
\end{tabular}
 
\caption{Expected values of $\alpha^\infty_{M,j}(8)$}
\label{tab8alphaInfMj}
\end{table}  

\begin{table}[h] \center

\begin{tabular}{||c||c|c|c|c||}
\hhline{|t:=:t:====:t|}
  \multicolumn{1}{||c||}{$M$}
&
  \multicolumn{1}{c|}{$j=1,\,7$}
&
  \multicolumn{1}{c|}{$j=2,\,6$}
&
  \multicolumn{1}{c|}{$j=3,\,5$}
&
  \multicolumn{1}{c||}{$j=4$}
\\
\hhline{||-||-|-|-|-||}

0
&$0$
&$8c_8^2$
&$4$
&$-16c_8^2+8$\\ 
1
&$6$
&$0$
&$2$
&$0$\\ 
2
&$2$
&$4c_8^2+4$
&$0$
&$-8c_8^2+4$\\ 
3
&$3$
&$4$
&$1$
&$0$\\ 
4
&$4$
&$8c_8^2$
&$0$
&$-16c_8^2+8$\\ 
5
&$2$
&$0$
&$6$
&$0$\\ 
6
&$0$
&$-\infty$
&$+\infty$
&$-\infty$\\ 
7
&$-\infty$
&$+\infty$
&$-\infty$
&$+\infty$\\ 
\hhline{|b:=:b:====:b|}
\end{tabular}
 
\caption{Expected values of $\sqrt{8}\beta^\infty_{M,j}(8)$}
\label{tab8betaInfMj}
\end{table}
 
\begin{table}[h] \center

\begin{tabular}{||c||c|c|c|c||}
\hhline{|t:=:t:====:t|}
  \multicolumn{1}{||c||}{$M$}
&
  \multicolumn{1}{c|}{$j=1,\,8$}
&
  \multicolumn{1}{c|}{$j=2,\,7$}
&
  \multicolumn{1}{c|}{$j=3,\,6$}
&
  \multicolumn{1}{c||}{$j=4,\,5$}
\\
\hhline{||-||-|-|-|-||}

0
&$0$
&$c_9+1$
&$-\frac{1}{2}$
&$-c_9+\frac{1}{2}$\\ 
1
&$-\frac{10c_9^2}{3}+\frac{c_9}{3}+\frac{7}{6}$
&$0$
&$\frac{1}{2}$
&$\frac{10c_9^2}{3}-\frac{c_9}{3}-\frac{2}{3}$\\ 
2
&$-\infty$
&$+\infty$
&$\frac{1}{2}$
&$-\infty$\\ 
3, 4
&$-\frac{2c_9^2}{3}-\frac{4c_9}{3}+\frac{4}{3}$
&$\frac{2c_9^2}{3}+\frac{4c_9}{3}-\frac{5}{6}$
&$\frac{1}{2}$
&$0$\\ 
5
&$+\infty$
&$-\infty$
&$\frac{1}{2}$
&$+\infty$\\ 
6
&$-2c_9^2+c_9+\frac{3}{2}$
&$0$
&$-\frac{1}{2}$
&$2c_9^2-c_9$\\ 
7
&$0$
&$+\infty$
&$-\infty$
&$+\infty$\\ 
8
&$+\infty$
&$-\infty$
&$+\infty$
&$-\infty$\\ 
\hhline{|b:=:b:====:b|}
\end{tabular}
 
\caption{Expected values of $\alpha^\infty_{M,j}(9)$}
\label{tab9alphaInfMj}
\end{table}  

\begin{table}[h] \center

\begin{tabular}{||c||c|c|c|c||}
\hhline{|t:=:t:====:t|}
  \multicolumn{1}{||c||}{$M$}
&
  \multicolumn{1}{c|}{$j=1,\,8$}
&
  \multicolumn{1}{c|}{$j=2,\,7$}
&
  \multicolumn{1}{c|}{$j=3,\,6$}
&
  \multicolumn{1}{c||}{$j=4,\,5$}
\\
\hhline{||-||-|-|-|-||}

0
&$0$
&$3c_9+3$
&$\frac{9}{2}$
&$-3c_9+\frac{3}{2}$\\ 
1
&$-2c_9^2+5c_9+\frac{11}{2}$
&$0$
&$\frac{3}{2}$
&$2c_9^2-5c_9+2$\\ 
2
&$-\infty$
&$+\infty$
&$\frac{3}{2}$
&$-\infty$\\ 
3, 4
&$-10c_9^2+4c_9+8$
&$10c_9^2-4c_9-\frac{1}{2}$
&$\frac{3}{2}$
&$0$\\ 
5
&$+\infty$
&$-\infty$
&$\frac{3}{2}$
&$+\infty$\\ 
6
&$-6c_9^2+3c_9+\frac{9}{2}$
&$0$
&$\frac{9}{2}$
&$6c_9^2-3c_9$\\ 
7
&$0$
&$-\infty$
&$+\infty$
&$-\infty$\\ 
8
&$-\infty$
&$+\infty$
&$-\infty$
&$+\infty$\\ 
\hhline{|b:=:b:====:b|}
\end{tabular}
 
\caption{Expected values of $\sqrt{9}\beta^\infty_{M,j}(9)$}
\label{tab9betaInfMj}
\end{table}

\begin{table}[h] \center

\begin{tabular}{||c||c|c|c|c|c||}
\hhline{|t:=:t:=====:t|}
  \multicolumn{1}{||c||}{$M$}
&
  \multicolumn{1}{c|}{$j=1,\,9$}
&
  \multicolumn{1}{c|}{$j=2,\,8$}
&
  \multicolumn{1}{c|}{$j=3,\,7$}
&
  \multicolumn{1}{c|}{$j=4,\,6$}
&
  \multicolumn{1}{c||}{$j=5$}
\\
\hhline{||-||-|-|-|-|-||}

0, 7
&$0$
&$0$
&$16c_{10}^2-4$
&$-48c_{10}^2+16$
&$64c_{10}^2-22$\\ 
1, 6
&$-8c_{10}^2+4$
&$0$
&$0$
&$8c_{10}^2-4$
&$2$\\ 
2, 5
&$16c_{10}^2-4$
&$-48c_{10}^2+16$
&$0$
&$0$
&$64c_{10}^2-22$\\ 
3, 4
&$-4c_{10}^2+3$
&$1$
&$4c_{10}^2-3$
&$0$
&$0$\\ 
8
&$0$
&$+\infty$
&$-\infty$
&$+\infty$
&$-\infty$\\ 
9
&$+\infty$
&$-\infty$
&$+\infty$
&$-\infty$
&$+\infty$\\ 
\hhline{|b:=:b:=====:b|}
\end{tabular}
\nocaption{Expected values of $\alpha^\infty_{M,j}(10)$}
 
\caption{Expected values of $\alpha^\infty_{M,j}(10)$}
\label{tab10alphaInfMj}
\end{table}  

\begin{table}[h] \center

\begin{tabular}{||c||c|c|c|c|c||}
\hhline{|t:=:t:=====:t|}
  \multicolumn{1}{||c||}{$M$}
&
  \multicolumn{1}{c|}{$j=1,\,9$}
&
  \multicolumn{1}{c|}{$j=2,\,8$}
&
  \multicolumn{1}{c|}{$j=3,\,7$}
&
  \multicolumn{1}{c|}{$j=4,\,6$}
&
  \multicolumn{1}{c||}{$j=5$}
\\
\hhline{||-||-|-|-|-|-||}

0, 7
&$0$
&$0$
&$64c_{10}^2-20$
&$-160c_{10}^2+60$
&$192c_{10}^2-60$\\ 
1, 6
&$16c_{10}^2$
&$0$
&$0$
&$0$
&$-32c_{10}^2+20$\\ 
2, 5
&$64c_{10}^2-20$
&$-160c_{10}^2+60$
&$0$
&$0$
&$192c_{10}^2-60$\\ 
3, 4
&$-8c_{10}^2+10$
&$5$
&$8c_{10}^2-5$
&$0$
&$0$\\ 
8
&$0$
&$+\infty$
&$-\infty$
&$+\infty$
&$-\infty$\\ 
9
&$+\infty$
&$-\infty$
&$+\infty$
&$-\infty$
&$+\infty$\\ 
\hhline{|b:=:b:=====:b|}
\end{tabular}
\nocaption{Expected values of $\sqrt{10}\beta^\infty_{M,j}(10)$}
 
\caption{Expected values of $\sqrt{10}\beta^\infty_{M,j}(10)$}
\label{tab10betaInfMj}
\end{table} 

\begin{table}[h] \center

\begin{tabular}{||c||c|c||}
\hhline{|t:=:t:==:t|}
  \multicolumn{1}{||c||}{$M$}
&
  \multicolumn{1}{c|}{$j=1,\,10$}
&
  \multicolumn{1}{c||}{$j=2,\,9$}
\\
\hhline{||-||-|-||}

0, 8
&$0$
&$0$\\ 
1, 7
&$-16c_{11}^4-8c_{11}^3+4c_{11}^2+2c_{11}+1$
&$0$\\ 
2, 6
&$-8c_{11}^3+8c_{11}^2+2c_{11}-1$
&$-8c_{11}^3+4c_{11}^2$\\ 
3, 4, 5
&$-32c_{11}^4-8c_{11}^3+28c_{11}^2+8c_{11}$
&$48c_{11}^4+16c_{11}^3-44c_{11}^2-16c_{11}+2$\\ 
9
&$0$
&$+\infty$\\ 
10
&$+\infty$
&$-\infty$\\ 
\hhline{|b:=:b:==:b|}
\end{tabular}
 
\\[2mm]

\begin{tabular}{||c||c|c||}
\hhline{|t:=:t:==:t|}
  \multicolumn{1}{||c||}{$M$}
&
  \multicolumn{1}{c|}{$j=3,\,8$}
&
  \multicolumn{1}{c||}{$j=4,\,7$}
\\
\hhline{||-||-|-||}

0, 8
&$16c_{11}^4+8c_{11}^3-12c_{11}^2-4c_{11}+2$
&$-16c_{11}^4-16c_{11}^3+4c_{11}^2+8c_{11}+2$\\ 
1, 7
&$0$
&$8c_{11}^3+4c_{11}^2-4c_{11}-1$\\ 
2, 6
&$0$
&$0$\\ 
3, 4, 5
&$-16c_{11}^4-8c_{11}^3+16c_{11}^2+8c_{11}-1$
&$0$\\ 
9
&$-\infty$
&$+\infty$\\ 
10
&$+\infty$
&$-\infty$\\ 
\hhline{|b:=:b:==:b|}
\end{tabular}
 
\label{tab11alphaInfMj}
\\[2mm]

\begin{tabular}{||c||c||}
\hhline{|t:=:t:=:t|}
  \multicolumn{1}{||c||}{$M$}
&
  \multicolumn{1}{c||}{$j=5,\,6$}
\\
\hhline{||-||-||}

0, 8
&$8c_{11}^3+8c_{11}^2-4c_{11}-3$\\ 
1, 7
&$16c_{11}^4-8c_{11}^2+2c_{11}+1$\\ 
2, 6
&$16c_{11}^3-12c_{11}^2-2c_{11}+2$\\ 
3, 4, 5
&$0$\\ 
9
&$-\infty$\\ 
10
&$+\infty$\\ 
\hhline{|b:=:b:=:b|}
\end{tabular}
\nocaption{Expected values of $\alpha^\infty_{M,j}(11)$}
 
\caption{Expected values of $\alpha^\infty_{M,j}(11)$}
\end{table}  

\begin{table}[h] \center

\begin{tabular}{||c||c|c||}
\hhline{|t:=:t:==:t|}
  \multicolumn{1}{||c||}{$M$}
&
  \multicolumn{1}{c|}{$j=1,\,10$}
&
  \multicolumn{1}{c||}{$j=2,\,9$}
\\
\hhline{||-||-|-||}

0, 8
&$0$
&$0$\\ 
1, 7
&$80c_{11}^4+40c_{11}^3-52c_{11}^2-14c_{11}+9$
&$0$\\ 
2, 6
&$32c_{11}^4-24c_{11}^3+2c_{11}+1$
&$-40c_{11}^3+12c_{11}^2+12c_{11}+4$\\ 
3, 4, 5
&$128c_{11}^4+24c_{11}^3-124c_{11}^2-28c_{11}+14$
&$-208c_{11}^4-16c_{11}^3+188c_{11}^2+32c_{11}-8$\\ 
9
&$0$
&$+\infty$\\ 
10
&$+\infty$
&$-\infty$\\ 
\hhline{|b:=:b:==:b|}
\end{tabular}
 
\\[2mm]

\begin{tabular}{||c||c|c||}
\hhline{|t:=:t:==:t|}
  \multicolumn{1}{||c||}{$M$}
&
  \multicolumn{1}{c|}{$j=3,\,8$}
&
  \multicolumn{1}{c||}{$j=4,\,7$}
\\
\hhline{||-||-|-||}

0, 8
&$80c_{11}^4+8c_{11}^3-60c_{11}^2+10$
&$-48c_{11}^4-48c_{11}^3+12c_{11}^2+20c_{11}+8$\\ 
1, 7
&$0$
&$-64c_{11}^4-8c_{11}^3+52c_{11}^2+4c_{11}-3$\\ 
2, 6
&$0$
&$0$\\ 
3, 4, 5
&$80c_{11}^4-8c_{11}^3-64c_{11}^2-4c_{11}+5$
&$0$\\ 
9
&$-\infty$
&$+\infty$\\ 
10
&$+\infty$
&$-\infty$\\ 
\hhline{|b:=:b:==:b|}
\end{tabular}

\\[2mm] 

\begin{tabular}{||c||c||}
\hhline{|t:=:t:=:t|}
  \multicolumn{1}{||c||}{$M$}
&
  \multicolumn{1}{c||}{$j=5,\,6$}
\\
\hhline{||-||-||}

0, 8
&$-32c_{11}^4+40c_{11}^3+48c_{11}^2-20c_{11}-7$\\ 
1, 7
&$-16c_{11}^4-32c_{11}^3+10c_{11}+5$\\ 
2, 6
&$-32c_{11}^4+64c_{11}^3-12c_{11}^2-14c_{11}+6$\\ 
3, 4, 5
&$0$\\ 
9
&$-\infty$\\ 
10
&$+\infty$\\ 
\hhline{|b:=:b:=:b|}
\end{tabular}
\nocaption{Expected values of $\sqrt{11}\beta^\infty_{M,j}(11)$}
 
\caption{Expected values of $\sqrt{11}\beta^\infty_{M,j}(11)$}
\label{tab11betaInfMj}
\end{table}

\clearpage
 
\begin{table}[h] \center

\begin{tabular}{||c||c|c|c|c|c|c||}
\hhline{|t:=:t:======:t|}
  \multicolumn{1}{||c||}{$M$}
&
  \multicolumn{1}{c|}{$j=1,\,11$}
&
  \multicolumn{1}{c|}{$j=2,\,10$}
&
  \multicolumn{1}{c|}{$j=3,\,9$}
&
  \multicolumn{1}{c|}{$j=4,\,8$}
&
  \multicolumn{1}{c|}{$j=5,\,7$}
&
  \multicolumn{1}{c||}{$j=6$}
\\
\hhline{||-||-|-|-|-|-|-||}

0
&$0$
&$0$
&$8c_{12}^2-1$
&$-12c_{12}^2+1$
&$4c_{12}^2+1$
&$0$\\ 
1
&$-\infty$
&$0$
&$0$
&$-2$
&$+\infty$
&$-6$\\ 
2
&$-2c_{12}^2+\frac{5}{2}$
&$-2$
&$0$
&$0$
&$2c_{12}^2-\frac{1}{2}$
&$2$\\ 
3
&$-\frac{2c_{12}^2}{3}-\frac{1}{6}$
&$\frac{2c_{12}^2}{3}+\frac{1}{6}$
&$-\frac{4c_{12}^2}{3}+\frac{13}{6}$
&$0$
&$0$
&$\frac{8c_{12}^2}{3}-\frac{7}{3}$\\ 
4, 5
&$-\frac{2c_{12}^2}{3}-\frac{1}{6}$
&$-c_{12}^2+\frac{7}{4}$
&$\frac{2c_{12}^2}{3}+\frac{1}{6}$
&$c_{12}^2-\frac{3}{4}$
&$0$
&$0$\\ 
6
&$-4c_{12}^2+3$
&$4c_{12}^2-3$
&$1$
&$0$
&$0$
&$0$\\ 
7
&$+\infty$
&$-2$
&$0$
&$0$
&$-\infty$
&$2$\\ 
8
&$-2c_{12}^2+\frac{5}{2}$
&$0$
&$0$
&$-2$
&$2c_{12}^2+\frac{7}{2}$
&$-6$\\ 
9
&$0$
&$0$
&$+\infty$
&$-\infty$
&$+\infty$
&$-\infty$\\ 
10
&$0$
&$+\infty$
&$-\infty$
&$+\infty$
&$-\infty$
&$+\infty$\\ 
11
&$+\infty$
&$-\infty$
&$+\infty$
&$-\infty$
&$+\infty$
&$-\infty$\\ 
\hhline{|b:=:b:======:b|}
\end{tabular}
\nocaption{Expected values of $\alpha^\infty_{M,j}(12)$}
 
\caption{Expected values of $\alpha^\infty_{M,j}(12)$}
\label{tab12alphaInfMj}
\end{table}  

\begin{table}[h] \center

\begin{tabular}{||c||c|c|c|c|c|c||}
\hhline{|t:=:t:======:t|}
  \multicolumn{1}{||c||}{$M$}
&
  \multicolumn{1}{c|}{$j=1,\,11$}
&
  \multicolumn{1}{c|}{$j=2,\,10$}
&
  \multicolumn{1}{c|}{$j=3,\,9$}
&
  \multicolumn{1}{c|}{$j=4,\,8$}
&
  \multicolumn{1}{c|}{$j=5,\,7$}
&
  \multicolumn{1}{c||}{$j=6$}
\\
\hhline{||-||-|-|-|-|-|-||}

0
&$0$
&$0$
&$24c_{12}^2$
&$-24c_{12}^2+6$
&$-24c_{12}^2+6$
&$48c_{12}^2$\\ 
1
&$+\infty$
&$0$
&$0$
&$12$
&$-\infty$
&$24$\\ 
2
&$-4c_{12}^2+11$
&$-4$
&$0$
&$0$
&$4c_{12}^2+1$
&$8$\\ 
3
&$4c_{12}^2-1$
&$-4c_{12}^2+9$
&$4c_{12}^2$
&$0$
&$0$
&$-8c_{12}^2+8$\\ 
4, 5
&$4c_{12}^2-1$
&$2c_{12}^2+\frac{7}{2}$
&$-4c_{12}^2+7$
&$-2c_{12}^2+\frac{5}{2}$
&$0$
&$0$\\ 
6
&$-8c_{12}^2+10$
&$8c_{12}^2-2$
&$8c_{12}^2-4$
&$0$
&$0$
&$-16c_{12}^2+16$\\ 
7
&$-\infty$
&$-4$
&$0$
&$0$
&$+\infty$
&$8$\\ 
8
&$12c_{12}^2-9$
&$0$
&$0$
&$12$
&$-12c_{12}^2-3$
&$24$\\ 
9
&$0$
&$0$
&$-\infty$
&$+\infty$
&$-\infty$
&$+\infty$\\ 
10
&$0$
&$-\infty$
&$+\infty$
&$-\infty$
&$+\infty$
&$-\infty$\\ 
11
&$-\infty$
&$+\infty$
&$-\infty$
&$+\infty$
&$-\infty$
&$+\infty$\\ 
\hhline{|b:=:b:======:b|}
\end{tabular}
\nocaption{Expected values of $\sqrt{12}\beta^\infty_{M,j}(12)$}
 
\caption{Expected values of $\sqrt{12}\beta^\infty_{M,j}(12)$}
\label{tab12betaInfMj}
\end{table}

  \clearpage
 
      \begin{table}[h] \center
 
\begin{tabular}{||r|r|r||c||r|r||c||}
\hhline{|t:===:t:=:t:==:t:=:t|}
  \multicolumn{1}{||c}{$q$}
&
  \multicolumn{1}{|c}{$M$}
&
  \multicolumn{1}{|c||}{$j$}
&
  \multicolumn{1}{c||}{$\alpha^\infty_{M,j}(q)$}
&
  \multicolumn{1}{c|}{$L$}
&
  \multicolumn{1}{c||}{$N$}
&
  \multicolumn{1}{c||}{$\Delta$}
\\
\hhline{||-|-|-||-||-|-||-||}$5$
&
$2$
&
$1$
&
$-c_{5}+1$
&
$15$
&
$77$
&
$-6.5062\dotsc\cdot 10^{-41}$\\ 
$5$
&
$2$
&
$1$
&
$-c_{5}+1$
&
$18$
&
$92$
&
$-5.0740\dotsc\cdot 10^{-49}$\\ 
$6$
&
$1$
&
$2$
&
$0$
&
$12$
&
$73$
&
$-1.7859\dotsc\cdot 10^{-31}$\\ 
$6$
&
$1$
&
$2$
&
$0$
&
$16$
&
$97$
&
$-2.9089\dotsc\cdot 10^{-42}$\\ 
$7$
&
$4$
&
$1$
&
$-4c_{7}^2+2c_{7}+1$
&
$8$
&
$60$
&
$-2.3692\dotsc\cdot 10^{-21}$\\ 
$7$
&
$4$
&
$1$
&
$-4c_{7}^2+2c_{7}+1$
&
$12$
&
$88$
&
$-4.2205\dotsc\cdot 10^{-32}$\\ 
$8$
&
$3$
&
$1$
&
$-c_{8}^2+\frac{1}{2}$
&
$8$
&
$67$
&
$-5.9491\dotsc\cdot 10^{-22}$\\ 
$8$
&
$3$
&
$1$
&
$-c_{8}^2+\frac{1}{2}$
&
$12$
&
$99$
&
$-1.0780\dotsc\cdot 10^{-32}$\\ 
$9$
&
$3$
&
$7$
&
$\frac{2c_{9}^2}{3}+\frac{4c_{9}}{3}-\frac{5}{6}$
&
$12$
&
$111$
&
$-6.6922\dotsc\cdot 10^{-34}$\\ 
$10$
&
$0$
&
$4$
&
$-48c_{10}^2+16$
&
$12$
&
$120$
&
$-7.9034\dotsc\cdot 10^{-29}$\\ 
$11$
&
$3$
&
$10$
&
$-32c_{11}^4-8c_{11}^3+28c_{11}^2+8c_{11}$
&
$11$
&
$124$
&
$-1.8409\dotsc\cdot 10^{-87}$\\ 
$12$
&
$3$
&
$6$
&
$\frac{8c_{12}^2}{3}-\frac{7}{3}$
&
$10$
&
$123$
&
$-7.6413\dotsc\cdot 10^{-27}$\\ 
\hhline{|b:===:b:=:b:==:b:=:b|}
\end{tabular}
\nocaption{Differenses $\Delta=\alpha_{N,0,j}(0,q)-\alpha^\infty_{M,j}(q)$ for  $N=qL+M$}
 
\caption{Differences $\Delta=\alpha_{N,0,j}(0,q)-\alpha^\infty_{M,j}(q)$ for  $N=qL+M$}
\label{tabalphalim}
\end{table}

      \begin{table}[h] \center
 
\begin{tabular}{||r|r|r||c||r|r||c||}
\hhline{|t:===:t:=:t:==:t:=:t|}
  \multicolumn{1}{||c}{$q$}
&
  \multicolumn{1}{|c}{$M$}
&
  \multicolumn{1}{|c||}{$j$}
&
  \multicolumn{1}{c||}{$\beta^\infty_{M,j}(q)$}
&
  \multicolumn{1}{c|}{$L$}
&
  \multicolumn{1}{c||}{$N$}
&
  \multicolumn{1}{c||}{$\Delta$}
\\
\hhline{||-|-|-||-||-|-||-||}$5$
&
$2$
&
$1$
&
$-c_{5}+4$
&
$15$
&
$77$
&
$1.4548\dotsc\cdot 10^{-40}$\\ 
$5$
&
$2$
&
$1$
&
$-c_{5}+4$
&
$18$
&
$92$
&
$1.1345\dotsc\cdot 10^{-48}$\\ 
$6$
&
$1$
&
$2$
&
$0$
&
$12$
&
$73$
&
$-3.5514\dotsc\cdot 10^{-31}$\\ 
$6$
&
$1$
&
$2$
&
$0$
&
$16$
&
$97$
&
$-5.7928\dotsc\cdot 10^{-42}$\\ 
$7$
&
$4$
&
$1$
&
$-4c_{7}^2+6c_{7}+3$
&
$8$
&
$60$
&
$6.8317\dotsc\cdot 10^{-21}$\\ 
$7$
&
$4$
&
$1$
&
$-4c_{7}^2+6c_{7}+3$
&
$12$
&
$88$
&
$1.2181\dotsc\cdot 10^{-31}$\\ 
$8$
&
$3$
&
$1$
&
$3$
&
$8$
&
$67$
&
$-2.4687\dotsc\cdot 10^{-22}$\\ 
$8$
&
$3$
&
$1$
&
$3$
&
$12$
&
$99$
&
$-4.4005\dotsc\cdot 10^{-33}$\\ 
$9$
&
$3$
&
$7$
&
$10c_{9}^2-4c_{9}-\frac{1}{2}$
&
$12$
&
$111$
&
$2.0076\dotsc\cdot 10^{-33}$\\ 
$10$
&
$0$
&
$4$
&
$-160c_{10}^2+60$
&
$12$
&
$120$
&
$-2.5032\dotsc\cdot 10^{-28}$\\ 
$11$
&
$3$
&
$10$
&
$128c_{11}^4+24c_{11}^3-124c_{11}^2-28c_{11}+14$
&
$11$
&
$124$
&
$-5.1136\dotsc\cdot 10^{-87}$\\ 
$12$
&
$3$
&
$6$
&
$-8c_{12}^2+8$
&
$10$
&
$123$
&
$-2.6470\dotsc\cdot 10^{-26}$\\ 
\hhline{|b:===:b:=:b:==:b:=:b|}
\end{tabular}
\nocaption{Differenses $\Delta=\sqrt{q}\beta_{N,0,j}(0,q)-\beta^\infty_{M,j}(q)$ for  $N=qL+M$}
 
\caption{Differences $\Delta=\sqrt{q}\beta_{N,0,j}(0,q)-\beta^\infty_{M,j}(q)$ for  $N=qL+M$}
\label{tabbetalim}
\end{table}

  \subsection{Other limiting values}
  
Examination of Tables  \ref{tab2alphaInfMj}--\ref{tab12betaInfMj}
shows that for $q=1\dots,12$ and $M=0,\dots,q-1$,
if all the values of
  $\alpha^\infty_{M,0}(q),\dots, \alpha^\infty_{M,q-1}(q)$
  are finite, then
  \begin{equation}\label{suma0}
    \sum_{j=0}^{q-1}\alpha^\infty_{M,j}(q)=0.
\end{equation}
Similarly, if all the values of
  $\beta^\infty_{M,0}(q),\dots, \beta^\infty_{M,q-1}(q)$
  are finite, then
  \begin{equation}\label{sumb2q}
    \sum_{j=0}^{q-1}\sqrt{q}\beta^\infty_{M,j}(q)=2q.
  \end{equation}
The numerical data suggest the following
more general surmise.
 \begin{myconjecture}  \label{conjsum1}
 For all $q$ and $L$, 
  \begin{equation}\label{sumAinf1}
    \sum_{j=0}^{q-1}\alpha_{N,0,L+j}(0,q)\rightarrow 0
\qquad
\text{and}\qquad 
    \sum_{j=0}^{q-1}\sqrt{q}\beta_{N,0,L+j}(0,q)\rightarrow 2q
  \end{equation}
as $N$ tends to infinity.  
\end{myconjecture}
Tables \ref{tabsuminfA1} and \ref{tabsuminfB1}
support Conjecture \ref{conjsum1}.

\begin{table}[h]\center
  \begin{tabular}{||r|r|r||c||l||}
  \hhline{|t:===:t:=:t:=:t|}
    
  \multicolumn{1}{||r}{$q$}
&
  \multicolumn{1}{|r}{$L$}
&
  \multicolumn{1}{|r||}{$N$}
&
  \multicolumn{1}{c||}{$  \alpha^\infty_{N,L}(q)$}
&
  \multicolumn{1}{c||}{$ \sum_{j=0}^{q-1}{ \alpha_{N,0,L+j}(q)}$}
\\

  \hhline{|:===::=::=:|}
    $8$
  &
    $2$
  &
    $15$
  &
    $-\infty$
  &
    $-7.18248\dotsc\cdot 10^{-20}$
  \\ 

    $8$
  &
    $2$
  &
    $23$
  &
    $-\infty$
  &
    $-2.85712\dotsc\cdot 10^{-62}$
  \\ 

    $8$
  &
    $2$
  &
    $31$
  &
    $-\infty$
  &
    $-2.42848\dotsc\cdot 10^{-126}$
  \\ 

    $8$
  &
    $2$
  &
    $39$
  &
    $-\infty$
  &
    $-3.70631\dotsc\cdot 10^{-212}$
  \\ 

    $9$
  &
    $2$
  &
    $25$
  &
    $+\infty$
  &
    $-1.39268\dotsc\cdot 10^{-68}$
  \\ 

    $9$
  &
    $2$
  &
    $34$
  &
    $+\infty$
  &
    $\phantom{-}2.43214\dotsc\cdot 10^{-140}$
  \\ 

    $9$
  &
    $2$
  &
    $43$
  &
    $+\infty$
  &
    $-1.21692\dotsc\cdot 10^{-236}$
  \\ 

    $9$
  &
    $2$
  &
    $52$
  &
    $+\infty$
  &
    $\phantom{-}1.72256\dotsc\cdot 10^{-357}$
  \\ 

  \hhline{|b:===:b:=:b:=:b|}
\end{tabular}
 
   \caption{Data supporting Conjecture \ref{conjsum1}}
  \label{tabsuminfA1}
\end{table}

\begin{table}[h]\center
   \begin{tabular}{||r|r|r||c||l||}
  \hhline{|t:===:t:=:t:=:t|}
    
  \multicolumn{1}{||r}{$q$}
&
  \multicolumn{1}{|r}{$L$}
&
  \multicolumn{1}{|r||}{$N$}
&
  \multicolumn{1}{c||}{$  \beta^\infty_{N,L}(q)$}
&
  \multicolumn{1}{c||}{$ \sum_{j=0}^{q-1}{ \sqrt{q}\beta_{N,0,L+j}(q)}$}
\\

  \hhline{|:===::=::=:|}
    $7$
  &
    $1$
  &
    $13$
  &
    $+\infty$
  &
    $14+3.08004\dotsc\cdot 10^{-18}$
  \\ 

    $7$
  &
    $1$
  &
    $20$
  &
    $+\infty$
  &
    $14-3.95278\dotsc\cdot 10^{-55}$
  \\ 

    $7$
  &
    $1$
  &
    $27$
  &
    $+\infty$
  &
    $14+4.49935\dotsc\cdot 10^{-111}$
  \\ 

    $7$
  &
    $1$
  &
    $34$
  &
    $+\infty$
  &
    $14-4.43640\dotsc\cdot 10^{-186}$
  \\ 

  \hhline{||-|-|-||-||-||}
    $9$
  &
    $1$
  &
    $20$
  &
    $-\infty$
  &
    $18-2.00702\dotsc\cdot 10^{-36}$
  \\ 

    $9$
  &
    $1$
  &
    $29$
  &
    $-\infty$
  &
    $18+5.50037\dotsc\cdot 10^{-92}$
  \\ 

    $9$
  &
    $1$
  &
    $38$
  &
    $-\infty$
  &
    $18-4.83907\dotsc\cdot 10^{-172}$
  \\ 

    $9$
  &
    $1$
  &
    $47$
  &
    $-\infty$
  &
    $18+1.28622\dotsc\cdot 10^{-276}$
  \\ 

  \hhline{|b:===:b:=:b:=:b|}
\end{tabular}
 
   \caption{Data supporting Conjecture \ref{conjsum1}}
  \label{tabsuminfB1}
\end{table}

According to the definitions \eqref{alims},
if  $\alpha_{M,j}^\infty(q)$ 
is different from zero and infinity, then 
\begin{eqnarray}\label{limgenA1}
 \frac{ \alpha_{qL+M,0,j}(0,q)}
  { \alpha_{M,j}^\infty(q)}&\rightarrow&1\label{limgenA}
\end{eqnarray}
as $L$ tends to infinity. Thus, according to Conjecture \ref{perij},
if $j_1\equiv j_2\bmod{q}$, then 
\begin{equation}\label{limgenA2}
  \frac{ \alpha_{qL+M,0,j_2}(0,q)}{ \alpha_{qL+M,0,j_1}(0,q)}\rightarrow 1. 
\end{equation}
The  numerical data (see Table \ref{zerozeroAeq}) suggest that the same is true 
even if 
the limiting values
are equal to zero or infinity.
Table \ref{zerozeroBeq} demonstrates similar phenomenon
for~$\beta_{N,0,j}(0,q)$.

\begin{myconjecture}  \label{conjlimgen1}
{For every $q$, $M,$  $j_1$ and $j_2$,  
if $ j_1\equiv  j_2\bmod{ q}
$ 
then} 
 \begin{equation}\label{limgen1}
  \frac{ \alpha_{qL+M,0,j_2}(0,q)}{ \alpha_{qL+M,0,j_1}(0,q)}\rightarrow 1
\qquad \text{and} \qquad 
  \frac{ \beta_{qL+M,0,j_2}(0,q)}{ \beta_{qL+M,0,j_1}(0,q)}\rightarrow
1
\end{equation}
as $L$ tends to infinity.
\end{myconjecture}

\begin{table}[h]\center
   \begin{tabular}{||r|r|r|r|r||c||l||}
  \hhline{|t:=====:t:=:t:=:t|}
    
  \multicolumn{1}{||r}{$q$}
&
  \multicolumn{1}{|r}{$M$}
&
  \multicolumn{1}{|r}{$j_1$}
&
  \multicolumn{1}{|r}{$j_2$}
&
  \multicolumn{1}{|r||}{$L$}
&
  \multicolumn{1}{c||}{$  \alpha^\infty_{M,j_1}(q)$}
&
  \multicolumn{1}{c||}{$ \frac{ \alpha_{qL+M,0,j_2}(q)}{ \alpha_{qL+M,0,j_1}(q)}-1$}
\\

  \hhline{|:=====::=::=:|}
    $7$
  &
    $0$
  &
    $1$
  &
    $8$
  &
    $5$
  &
    $0$
  &
    $\phantom{-}3.93820\dotsc\cdot 10^{-198}$
  \\ 

    $7$
  &
    $0$
  &
    $1$
  &
    $8$
  &
    $7$
  &
    $0$
  &
    $\phantom{-}8.51409\dotsc\cdot 10^{-411}$
  \\ 

    $7$
  &
    $0$
  &
    $1$
  &
    $8$
  &
    $9$
  &
    $0$
  &
    $\phantom{-}7.81033\dotsc\cdot 10^{-700}$
  \\ 

    $7$
  &
    $0$
  &
    $1$
  &
    $8$
  &
    $11$
  &
    $0$
  &
    $\phantom{-}2.99871\dotsc\cdot 10^{-1065}$
  \\ 

    $7$
  &
    $0$
  &
    $1$
  &
    $8$
  &
    $13$
  &
    $0$
  &
    $\phantom{-}4.73991\dotsc\cdot 10^{-1507}$
  \\ 

  \hhline{||-|-|-|-|-||-||-||}
    $9$
  &
    $1$
  &
    $2$
  &
    $11$
  &
    $5$
  &
    $0$
  &
    $-1.32200\dotsc\cdot 10^{-260}$
  \\ 

    $9$
  &
    $1$
  &
    $2$
  &
    $11$
  &
    $7$
  &
    $0$
  &
    $-6.56976\dotsc\cdot 10^{-538}$
  \\ 

    $9$
  &
    $1$
  &
    $2$
  &
    $11$
  &
    $9$
  &
    $0$
  &
    $-2.15423\dotsc\cdot 10^{-913}$
  \\ 

    $9$
  &
    $1$
  &
    $2$
  &
    $11$
  &
    $11$
  &
    $0$
  &
    $-4.43429\dotsc\cdot 10^{-1387}$
  \\ 

  \hhline{||-|-|-|-|-||-||-||}
    $8$
  &
    $6$
  &
    $2$
  &
    $10$
  &
    $5$
  &
    $+\infty$
  &
    $\phantom{-}3.35957\dotsc\cdot 10^{-345}$
  \\ 

    $8$
  &
    $6$
  &
    $2$
  &
    $10$
  &
    $7$
  &
    $+\infty$
  &
    $\phantom{-}2.65729\dotsc\cdot 10^{-636}$
  \\ 

    $8$
  &
    $6$
  &
    $2$
  &
    $10$
  &
    $9$
  &
    $+\infty$
  &
    $\phantom{-}9.67441\dotsc\cdot 10^{-1015}$
  \\ 

  \hhline{||-|-|-|-|-||-||-||}
    $10$
  &
    $8$
  &
    $3$
  &
    $13$
  &
    $3$
  &
    $-\infty$
  &
    $\phantom{-}5.19665\dotsc\cdot 10^{-172}$
  \\ 

    $10$
  &
    $8$
  &
    $3$
  &
    $13$
  &
    $5$
  &
    $-\infty$
  &
    $\phantom{-}3.54600\dotsc\cdot 10^{-424}$
  \\ 

    $10$
  &
    $8$
  &
    $3$
  &
    $13$
  &
    $7$
  &
    $-\infty$
  &
    $\phantom{-}1.79924\dotsc\cdot 10^{-785}$
  \\ 

    $10$
  &
    $8$
  &
    $3$
  &
    $13$
  &
    $9$
  &
    $-\infty$
  &
    $\phantom{-}6.65762\dotsc\cdot 10^{-1256}$
  \\ 

  \hhline{|b:=====:b:=:b:=:b|}
\end{tabular}
 
  \caption{Data supporting Conjecture \ref{conjlimgen1}}
  \label{zerozeroAeq}
\end{table}

\begin{table}\center
   \begin{tabular}{||r|r|r|r|r||c||l||}
  \hhline{|t:=====:t:=:t:=:t|}
    
  \multicolumn{1}{||r}{$q$}
&
  \multicolumn{1}{|r}{$M$}
&
  \multicolumn{1}{|r}{$j_1$}
&
  \multicolumn{1}{|r}{$j_2$}
&
  \multicolumn{1}{|r||}{$L$}
&
  \multicolumn{1}{c||}{$  \beta^\infty_{M,j_1}(q)$}
&
  \multicolumn{1}{c||}{$ \frac{ \beta_{qL+M,0,j_2}(q)}{ \beta_{qL+M,0,j_1}(q)}-1$}
\\

  \hhline{|:=====::=::=:|}
    $6$
  &
    $1$
  &
    $2$
  &
    $8$
  &
    $5$
  &
    $0$
  &
    $\phantom{-}2.33277\dotsc\cdot 10^{-183}$
  \\ 

    $6$
  &
    $1$
  &
    $2$
  &
    $8$
  &
    $7$
  &
    $0$
  &
    $\phantom{-}5.21090\dotsc\cdot 10^{-372}$
  \\ 

    $6$
  &
    $1$
  &
    $2$
  &
    $8$
  &
    $9$
  &
    $0$
  &
    $\phantom{-}4.14696\dotsc\cdot 10^{-626}$
  \\ 

    $6$
  &
    $1$
  &
    $2$
  &
    $8$
  &
    $11$
  &
    $0$
  &
    $\phantom{-}1.17306\dotsc\cdot 10^{-945}$
  \\ 

    $6$
  &
    $1$
  &
    $2$
  &
    $8$
  &
    $13$
  &
    $0$
  &
    $\phantom{-}1.15307\dotsc\cdot 10^{-1330}$
  \\ 

  \hhline{||-|-|-|-|-||-||-||}
    $7$
  &
    $1$
  &
    $2$
  &
    $9$
  &
    $5$
  &
    $0$
  &
    $\phantom{-}3.66642\dotsc\cdot 10^{-208}$
  \\ 

    $7$
  &
    $1$
  &
    $2$
  &
    $9$
  &
    $7$
  &
    $0$
  &
    $\phantom{-}2.63024\dotsc\cdot 10^{-426}$
  \\ 

    $7$
  &
    $1$
  &
    $2$
  &
    $9$
  &
    $9$
  &
    $0$
  &
    $\phantom{-}8.27575\dotsc\cdot 10^{-721}$
  \\ 

    $7$
  &
    $1$
  &
    $2$
  &
    $9$
  &
    $11$
  &
    $0$
  &
    $\phantom{-}1.09991\dotsc\cdot 10^{-1091}$
  \\ 

  \hhline{||-|-|-|-|-||-||-||}
    $7$
  &
    $6$
  &
    $1$
  &
    $8$
  &
    $5$
  &
    $+\infty$
  &
    $-1.10239\dotsc\cdot 10^{-291}$
  \\ 

    $7$
  &
    $6$
  &
    $1$
  &
    $8$
  &
    $7$
  &
    $+\infty$
  &
    $-2.70991\dotsc\cdot 10^{-542}$
  \\ 

    $7$
  &
    $6$
  &
    $1$
  &
    $8$
  &
    $9$
  &
    $+\infty$
  &
    $-2.51149\dotsc\cdot 10^{-869}$
  \\ 

    $7$
  &
    $6$
  &
    $1$
  &
    $8$
  &
    $11$
  &
    $+\infty$
  &
    $-8.93391\dotsc\cdot 10^{-1273}$
  \\ 

  \hhline{||-|-|-|-|-||-||-||}
    $9$
  &
    $2$
  &
    $1$
  &
    $10$
  &
    $3$
  &
    $-\infty$
  &
    $\phantom{-}3.65528\dotsc\cdot 10^{-97}$
  \\ 

    $9$
  &
    $2$
  &
    $1$
  &
    $10$
  &
    $5$
  &
    $-\infty$
  &
    $\phantom{-}4.90994\dotsc\cdot 10^{-287}$
  \\ 

    $9$
  &
    $2$
  &
    $1$
  &
    $10$
  &
    $7$
  &
    $-\infty$
  &
    $\phantom{-}4.14831\dotsc\cdot 10^{-575}$
  \\ 

  \hhline{|b:=====:b:=:b:=:b|}
\end{tabular}
 
   \caption{Data supporting Conjecture \ref{conjlimgen1}}
  \label{zerozeroBeq}
\end{table}

Similarly,  Conjecture \ref{consymm} implies 
\eqref{limgenA2}  whenever 
 $j_1+j_2 \equiv 0\bmod{q}$,
and the numerical data suggest (see Tables \ref{zerozeroAsum}) that  
 $\alpha_{M,j_1}^\infty(q)$ and $ \alpha_{M,j_2}^\infty(q)$
can be equal to zero or infinity as well.
Table \ref{zerozeroBsum} demonstrates similar phenomenon
for~$\beta_{N,0,j}(0,q)$.

\begin{myconjecture}  \label{conjlimgen2}
{For every $q$, $M,$  $j_1$ and $j_2$,  
if 
$ j_1+ j_2\equiv 0\bmod{ q}
$ 
then} 
 \begin{equation}\label{limgen2}
  \frac{ \alpha_{qL+M,0,j_2}(0,q)}{ \alpha_{qL+M,0,j_1}(0,q)}\rightarrow
1
\qquad \text{and} \qquad 
  \frac{ \beta_{qL+M,0,j_2}(0,q)}{ \beta_{qL+M,0,j_1}(0,q)}\rightarrow
1
\end{equation}
as $L$ tends to infinity.
\end{myconjecture}

\begin{table}[h]\center
   \begin{tabular}{||r|r|r|r|r||c||l||}
  \hhline{|t:=====:t:=:t:=:t|}
    
  \multicolumn{1}{||r}{$q$}
&
  \multicolumn{1}{|r}{$M$}
&
  \multicolumn{1}{|r}{$j_1$}
&
  \multicolumn{1}{|r}{$j_2$}
&
  \multicolumn{1}{|r||}{$N$}
&
  \multicolumn{1}{c||}{$  \alpha^\infty_{M,j_1}(q)$}
&
  \multicolumn{1}{c||}{$ \frac{ \alpha_{qN+M,0,j_2}(q)}{ \alpha_{qN+M,0,j_1}(q)}-1$}
\\

  \hhline{|:=====::=::=:|}
    $7$
  &
    $0$
  &
    $1$
  &
    $6$
  &
    $5$
  &
    $0$
  &
    $\phantom{-}8.39093\dotsc\cdot 10^{-201}$
  \\ 

    $7$
  &
    $0$
  &
    $1$
  &
    $6$
  &
    $7$
  &
    $0$
  &
    $\phantom{-}1.68951\dotsc\cdot 10^{-413}$
  \\ 

    $7$
  &
    $0$
  &
    $1$
  &
    $6$
  &
    $9$
  &
    $0$
  &
    $\phantom{-}1.51061\dotsc\cdot 10^{-702}$
  \\ 

    $7$
  &
    $0$
  &
    $1$
  &
    $6$
  &
    $11$
  &
    $0$
  &
    $\phantom{-}5.72953\dotsc\cdot 10^{-1068}$
  \\ 

    $7$
  &
    $0$
  &
    $1$
  &
    $6$
  &
    $13$
  &
    $0$
  &
    $\phantom{-}8.99531\dotsc\cdot 10^{-1510}$
  \\ 

  \hhline{||-|-|-|-|-||-||-||}
    $9$
  &
    $1$
  &
    $2$
  &
    $7$
  &
    $5$
  &
    $0$
  &
    $-4.64608\dotsc\cdot 10^{-266}$
  \\ 

    $9$
  &
    $1$
  &
    $2$
  &
    $7$
  &
    $7$
  &
    $0$
  &
    $-2.29737\dotsc\cdot 10^{-543}$
  \\ 

    $9$
  &
    $1$
  &
    $2$
  &
    $7$
  &
    $9$
  &
    $0$
  &
    $-7.52208\dotsc\cdot 10^{-919}$
  \\ 

    $9$
  &
    $1$
  &
    $2$
  &
    $7$
  &
    $11$
  &
    $0$
  &
    $-1.54744\dotsc\cdot 10^{-1392}$
  \\ 

  \hhline{||-|-|-|-|-||-||-||}
    $8$
  &
    $6$
  &
    $2$
  &
    $6$
  &
    $5$
  &
    $+\infty$
  &
    $\phantom{-}3.91745\dotsc\cdot 10^{-351}$
  \\ 

    $8$
  &
    $6$
  &
    $2$
  &
    $6$
  &
    $7$
  &
    $+\infty$
  &
    $\phantom{-}3.09265\dotsc\cdot 10^{-642}$
  \\ 

    $8$
  &
    $6$
  &
    $2$
  &
    $6$
  &
    $9$
  &
    $+\infty$
  &
    $\phantom{-}1.12525\dotsc\cdot 10^{-1020}$
  \\ 

  \hhline{||-|-|-|-|-||-||-||}
    $10$
  &
    $8$
  &
    $3$
  &
    $7$
  &
    $3$
  &
    $-\infty$
  &
    $\phantom{-}8.32708\dotsc\cdot 10^{-181}$
  \\ 

    $10$
  &
    $8$
  &
    $3$
  &
    $7$
  &
    $5$
  &
    $-\infty$
  &
    $\phantom{-}5.70848\dotsc\cdot 10^{-433}$
  \\ 

    $10$
  &
    $8$
  &
    $3$
  &
    $7$
  &
    $7$
  &
    $-\infty$
  &
    $\phantom{-}2.90808\dotsc\cdot 10^{-794}$
  \\ 

    $10$
  &
    $8$
  &
    $3$
  &
    $7$
  &
    $9$
  &
    $-\infty$
  &
    $\phantom{-}1.07854\dotsc\cdot 10^{-1264}$
  \\ 

  \hhline{|b:=====:b:=:b:=:b|}
\end{tabular}
 
  \caption{Data supporting Conjecture \ref{conjlimgen2}.}
  \label{zerozeroAsum}
\end{table}

\begin{table}[h]\center
   \begin{tabular}{||r|r|r|r|r||c||l||}
  \hhline{|t:=====:t:=:t:=:t|}
    
  \multicolumn{1}{||r}{$q$}
&
  \multicolumn{1}{|r}{$M$}
&
  \multicolumn{1}{|r}{$j_1$}
&
  \multicolumn{1}{|r}{$j_2$}
&
  \multicolumn{1}{|r||}{$N$}
&
  \multicolumn{1}{c||}{$  \beta^\infty_{M,j_1}(q)$}
&
  \multicolumn{1}{c||}{$ \frac{ \beta_{qN+M,0,j_2}(q)}{ \beta_{qN+M,0,j_1}(q)}-1$}
\\

  \hhline{|:=====::=::=:|}
    $6$
  &
    $1$
  &
    $2$
  &
    $4$
  &
    $5$
  &
    $0$
  &
    $\phantom{-}2.27340\dotsc\cdot 10^{-189}$
  \\ 

    $6$
  &
    $1$
  &
    $2$
  &
    $4$
  &
    $7$
  &
    $0$
  &
    $\phantom{-}4.23522\dotsc\cdot 10^{-378}$
  \\ 

    $6$
  &
    $1$
  &
    $2$
  &
    $4$
  &
    $9$
  &
    $0$
  &
    $\phantom{-}3.16022\dotsc\cdot 10^{-632}$
  \\ 

    $6$
  &
    $1$
  &
    $2$
  &
    $4$
  &
    $11$
  &
    $0$
  &
    $\phantom{-}8.66986\dotsc\cdot 10^{-952}$
  \\ 

    $6$
  &
    $1$
  &
    $2$
  &
    $4$
  &
    $13$
  &
    $0$
  &
    $\phantom{-}8.37832\dotsc\cdot 10^{-1337}$
  \\ 

  \hhline{||-|-|-|-|-||-||-||}
    $7$
  &
    $1$
  &
    $2$
  &
    $5$
  &
    $5$
  &
    $0$
  &
    $\phantom{-}1.38517\dotsc\cdot 10^{-213}$
  \\ 

    $7$
  &
    $1$
  &
    $2$
  &
    $5$
  &
    $7$
  &
    $0$
  &
    $\phantom{-}9.53090\dotsc\cdot 10^{-432}$
  \\ 

    $7$
  &
    $1$
  &
    $2$
  &
    $5$
  &
    $9$
  &
    $0$
  &
    $\phantom{-}2.95147\dotsc\cdot 10^{-726}$
  \\ 

    $7$
  &
    $1$
  &
    $2$
  &
    $5$
  &
    $11$
  &
    $0$
  &
    $\phantom{-}3.89260\dotsc\cdot 10^{-1097}$
  \\ 

  \hhline{||-|-|-|-|-||-||-||}
    $7$
  &
    $6$
  &
    $1$
  &
    $6$
  &
    $5$
  &
    $+\infty$
  &
    $-2.08349\dotsc\cdot 10^{-294}$
  \\ 

    $7$
  &
    $6$
  &
    $1$
  &
    $6$
  &
    $7$
  &
    $+\infty$
  &
    $-5.09156\dotsc\cdot 10^{-545}$
  \\ 

    $7$
  &
    $6$
  &
    $1$
  &
    $6$
  &
    $9$
  &
    $+\infty$
  &
    $-4.70732\dotsc\cdot 10^{-872}$
  \\ 

    $7$
  &
    $6$
  &
    $1$
  &
    $6$
  &
    $11$
  &
    $+\infty$
  &
    $-1.67244\dotsc\cdot 10^{-1275}$
  \\ 

  \hhline{||-|-|-|-|-||-||-||}
    $9$
  &
    $2$
  &
    $1$
  &
    $8$
  &
    $3$
  &
    $-\infty$
  &
    $\phantom{-}6.93014\dotsc\cdot 10^{-100}$
  \\ 

    $9$
  &
    $2$
  &
    $1$
  &
    $8$
  &
    $5$
  &
    $-\infty$
  &
    $\phantom{-}9.19843\dotsc\cdot 10^{-290}$
  \\ 

    $9$
  &
    $2$
  &
    $1$
  &
    $8$
  &
    $7$
  &
    $-\infty$
  &
    $\phantom{-}7.75569\dotsc\cdot 10^{-578}$
  \\ 

  \hhline{|b:=====:b:=:b:=:b|}
\end{tabular}
 
   \caption{Data supporting Conjecture \ref{conjlimgen2}.}
  \label{zerozeroBsum}
\end{table}

It is interesting to observe that the convergence in 
\eqref{limgen1} an \eqref{limgen2}  appears to be much faster than
then the convergence in \eqref{limgenA}.
 Tables \ref{fastcovA} and \ref{fastcovB} 
illustrate this phenomenon and 
its counterpart 
for   $\beta_{N,0,j}(0,q)$. 
\begin{table}[h]\center
   \begin{tabular}{||r|r|r|r|r||l||l||}
  \hhline{|t:=====:t:=:t:=:t|}
    
  \multicolumn{1}{||r}{$q$}
&
  \multicolumn{1}{|r}{$M$}
&
  \multicolumn{1}{|r}{$j_1$}
&
  \multicolumn{1}{|r}{$j_2$}
&
  \multicolumn{1}{|r||}{$L$}
&
  \multicolumn{1}{c||}{$ \frac{ \alpha_{qL+M,0,j_1}(q)}{ \alpha^\infty_{M,j_1}(q)}-1$}
&
  \multicolumn{1}{c||}{$ \frac{ \alpha_{qL+M,0,j_2}(q)}{ \alpha_{qL+M,0,j_1}(q)}-1$}
\\

  \hhline{|:=====::=::=:|}
    $4$
  &
    $2$
  &
    $1$
  &
    $3$
  &
    $5$
  &
    $\phantom{-}1.82489\dotsc\cdot 10^{-14}$
  &
    $-3.42997\dotsc\cdot 10^{-148}$
  \\ 

    $4$
  &
    $1$
  &
    $1$
  &
    $3$
  &
    $10$
  &
    $-1.92505\dotsc\cdot 10^{-27}$
  &
    $\phantom{-}1.92335\dotsc\cdot 10^{-539}$
  \\ 

    $5$
  &
    $1$
  &
    $1$
  &
    $4$
  &
    $8$
  &
    $-2.74429\dotsc\cdot 10^{-42}$
  &
    $-2.42904\dotsc\cdot 10^{-417}$
  \\ 

    $5$
  &
    $1$
  &
    $1$
  &
    $4$
  &
    $10$
  &
    $-5.23088\dotsc\cdot 10^{-53}$
  &
    $-4.85932\dotsc\cdot 10^{-658}$
  \\ 

    $6$
  &
    $2$
  &
    $2$
  &
    $4$
  &
    $8$
  &
    $\phantom{-}1.30423\dotsc\cdot 10^{-21}$
  &
    $\phantom{-}3.74313\dotsc\cdot 10^{-539}$
  \\ 

    $7$
  &
    $3$
  &
    $1$
  &
    $6$
  &
    $8$
  &
    $\phantom{-}3.90028\dotsc\cdot 10^{-21}$
  &
    $\phantom{-}2.59440\dotsc\cdot 10^{-632}$
  \\ 

    $7$
  &
    $3$
  &
    $1$
  &
    $6$
  &
    $12$
  &
    $\phantom{-}7.04072\dotsc\cdot 10^{-32}$
  &
    $\phantom{-}1.24261\dotsc\cdot 10^{-1406}$
  \\ 

    $10$
  &
    $4$
  &
    $2$
  &
    $8$
  &
    $8$
  &
    $\phantom{-}1.19614\dotsc\cdot 10^{-21}$
  &
    $-1.67492\dotsc\cdot 10^{-902}$
  \\ 

    $10$
  &
    $4$
  &
    $1$
  &
    $9$
  &
    $10$
  &
    $\phantom{-}5.24047\dotsc\cdot 10^{-27}$
  &
    $-4.97789\dotsc\cdot 10^{-1399}$
  \\ 

    $12$
  &
    $2$
  &
    $1$
  &
    $11$
  &
    $8$
  &
    $\phantom{-}2.60688\dotsc\cdot 10^{-20}$
  &
    $-6.95177\dotsc\cdot 10^{-1005}$
  \\ 

  \hhline{|b:=====:b:=:b:=:b|}
\end{tabular}
 
   \caption{Comparison  of the convergence rate in
  \eqref{limgenA} and  \eqref{limgen2}}
  \label{fastcovA}
\end{table}

\begin{table}[h]\center
  \begin{tabular}{||r|r|r|r|r||l||l||}
  \hhline{|t:=====:t:=:t:=:t|}
    
  \multicolumn{1}{||r}{$q$}
&
  \multicolumn{1}{|r}{$M$}
&
  \multicolumn{1}{|r}{$j_1$}
&
  \multicolumn{1}{|r}{$j_2$}
&
  \multicolumn{1}{|r||}{$L$}
&
  \multicolumn{1}{c||}{$ \frac{ \sqrt{q}\beta_{qL+M,0,j_1}(q)}{ \beta^\infty_{M,j_1}(q)}-1$}
&
  \multicolumn{1}{c||}{$ \frac{ \beta_{qL+M,0,j_2}(q)}{ \beta_{qL+M,0,j_1}(q)}-1$}
\\

  \hhline{|:=====::=::=:|}
    $4$
  &
    $2$
  &
    $1$
  &
    $3$
  &
    $5$
  &
    $-1.82489\dotsc\cdot 10^{-14}$
  &
    $-3.42997\dotsc\cdot 10^{-148}$
  \\ 

    $4$
  &
    $1$
  &
    $1$
  &
    $3$
  &
    $10$
  &
    $\phantom{-}6.41686\dotsc\cdot 10^{-28}$
  &
    $\phantom{-}6.41119\dotsc\cdot 10^{-540}$
  \\ 

    $5$
  &
    $1$
  &
    $1$
  &
    $4$
  &
    $8$
  &
    $\phantom{-}3.67271\dotsc\cdot 10^{-43}$
  &
    $-3.25080\dotsc\cdot 10^{-418}$
  \\ 

    $5$
  &
    $1$
  &
    $1$
  &
    $4$
  &
    $10$
  &
    $\phantom{-}7.00051\dotsc\cdot 10^{-54}$
  &
    $-6.50326\dotsc\cdot 10^{-659}$
  \\ 

    $6$
  &
    $2$
  &
    $2$
  &
    $4$
  &
    $8$
  &
    $\phantom{-}8.62134\dotsc\cdot 10^{-22}$
  &
    $-3.11917\dotsc\cdot 10^{-539}$
  \\ 

    $7$
  &
    $3$
  &
    $1$
  &
    $6$
  &
    $8$
  &
    $-2.73039\dotsc\cdot 10^{-22}$
  &
    $\phantom{-}6.14564\dotsc\cdot 10^{-633}$
  \\ 

    $7$
  &
    $3$
  &
    $1$
  &
    $6$
  &
    $12$
  &
    $-5.01133\dotsc\cdot 10^{-33}$
  &
    $\phantom{-}2.88398\dotsc\cdot 10^{-1407}$
  \\ 

    $10$
  &
    $4$
  &
    $2$
  &
    $8$
  &
    $8$
  &
    $\phantom{-}8.68413\dotsc\cdot 10^{-22}$
  &
    $\phantom{-}1.11498\dotsc\cdot 10^{-902}$
  \\ 

    $10$
  &
    $4$
  &
    $1$
  &
    $9$
  &
    $10$
  &
    $-1.59000\dotsc\cdot 10^{-27}$
  &
    $-3.71255\dotsc\cdot 10^{-1399}$
  \\ 

    $12$
  &
    $2$
  &
    $1$
  &
    $11$
  &
    $8$
  &
    $\phantom{-}8.99764\dotsc\cdot 10^{-21}$
  &
    $-2.10061\dotsc\cdot 10^{-1005}$
  \\ 

  \hhline{|b:=====:b:=:b:=:b|}
\end{tabular}
 
  \caption{The counterpart of Table \ref{fastcovA} for
   $\beta_{qL+M,0,j}(0,q)$.}
  \label{fastcovB}
\end{table}

\section*{Details of computation}
It is shown in Appendix C that the determinants of matrices \eqref{defB} 
do not vanish. However these matrices almost
degenerate, and calculation of $\alpha_{N,i,j}(d,q)$ and $\beta_{N,i,j}(d,q)$ should be performed with multiprecision.  
Actually this was done by a 
program  written in {\sc \mbox{Julia}}
programming language~\cite{JULIA};
package {\sc Nemo} \cite{Nemo}
was used for calculations in ball arithmetic
with guaranteed high precision.

Approximate numerical values of 
$\alpha_{M,j}^\infty(q)$ 
were calculate as common initial decimal digits 
of
$\alpha_{qL+M,0,j}(0,q)$ and $\alpha_{q(L+1)+M,0,j}(0,q)$
for sufficiently large $L$. Exact values of 
$\alpha_{M,j}^\infty(q)$  were guessed by LLL-algorithm
implemented in {\sc Mathematica} \cite{Mathematica}.
Similar technique was used for $\beta_{M,j}^\infty(q)$.

\section*{Summary}

We have introduced certain matrices $B_N(d,q)$ 
of size $2N\times 2N$; 
the entries to these matrices were suggested by a special
form of the functional equation for a Dirichlet $L$-function $L_\chi(s)$, $\chi$ being  a Dirichlet character modulo~$q$.

We have proven that these matrices are not singular. Numerical
calculations indicate that   matrices the  $B_N(0,q)$ and their 
inverses $B_N^{-1}(0,q)$ exhibits the following 
unexpected properties.

The formal definition of  $B_N(d,q)$ does not imply that
integer values of $q$ could play any special role. However, the
determinants of
  $B_N(0,q)$ vanish when $q$ is close to $4$, $5$, ... .
  
The entries in bottom row of matrices  $B_N^{-1}(0,q)$ are antisymmetrical, this implies the equality of the corresponding 
minors of the  matrices  $B_N(0,q)$. 

In addition to these \emph{exact} equalities between
the  entries of matrices 
  $B_N^{-1}(0,q)$, there are also many \emph{approximate} equalities.
In the left-hand half of the top rows of these matrices,  the
entries at distance $q$ are very close to each other; 
the same is true for the right-hand half of the same top rows.
Also there are similar approximate equalities between 
entries (at similar positions)~to
the top rows of matrices
 $B_{N_1}^{-1}(0,q)$  and  $B_{N_2}^{-1}(0,q)$ 
 provided that $N_1\equiv N_2\bmod{q}$. Moreover,
 there are (finite or infinite) limits of the values 
 of such entries as $N$ tends to $\infty,$ running along an
 arithmetic progression modulo~$q$. If the limits are finite, they belong to the field~$\mathbb{Q}(\cos(\pi/q))$ (up to a scaling factor of~$\sqrt{q}$).

\appendix
 \setcounter{section}{1}
\setcounter{equation}{0}
\numberwithin{equation}{section}

\clearpage
 \section*{Appendix A. Proof of equality \eqref{defEsum}}
 \label{appA}

 We need to show that \eqref{deflambda} is equal to
 \eqref{defElambda} with $E_m(x)$ defined by \eqref{defEsum}.
 To this end we prove more general equality
 \begin{eqnarray}\label{genidL}  
  \frac{\myd^m}{\myd \tau^m}  
  \mye^{-\frac{x}{\tau}}
  &=&   
  \tau^{-m} \mye^{-\frac{x}{\tau}}
    E_m\mleft(\frac{x}{\tau}\mright).
      \end{eqnarray}

The proof is by induction on $m$. Cases $m=0$ and $m=1$ are trivial.
Further  we have:
\begin{eqnarray}\label{genidLind}  
 \lefteqn{ \frac{\myd^{m+1}}{\myd \tau^{m+1}}  
  \mye^{-\frac{x}{\tau}}
  =}
  \\
  &=&
   \left(\frac{\myd^m}{\myd \tau^m}  
  \mye^{-\frac{x}{\tau}}\right)'
   \\
  &=&     
     \left(\tau^{-m} \mye^{-\frac{x}{\tau}}
    E_m\mleft(\frac{x}{\tau}\mright)\right)'\\
       &=&   
      \left(   \mye^{-\frac{x}{\tau}}
   \sum_{k=1}^m (-1)^{m+k}
  \frac{m!}{k!}\binom{m-1}{k-1}
  {x}^k{\tau}^{-m-k}\right)'
 \\
        &=&   
      {x}{\tau}^{-2}   \mye^{-\frac{x}{\tau}}
   \sum_{k=1}^m (-1)^{m+k}
  \frac{m!}{k!}\binom{m-1}{k-1}
  {x}^k{\tau}^{-m-k}
  \\&& 
 \quad +     \mye^{-\frac{x}{\tau}}
   \sum_{k=1}^m (-1)^{m+k+1}(m+k)
  \frac{m!}{k!}\binom{m-1}{k-1}
  {x}^k{\tau}^{-m-k-1}
  \\
        &=&   
      {x}{\tau}^{-2}   \mye^{-\frac{x}{\tau}}
   \sum_{k=2}^{m+1} (-1)^{m+k-1}
  \frac{m!}{(k-1)!}\binom{m-1}{k-2}
  {x}^{k-1}{\tau}^{-m-k+1}
  \\&& 
 \quad +     \mye^{-\frac{x}{\tau}}
   \sum_{k=1}^{m+1} (-1)^{m+k+1}(m+k)
  \frac{m!}{k!}\binom{m-1}{k-1}
  {x}^k{\tau}^{-m-k-1}
  \\
    &=&   %
  \tau^{-m-1} \mye^{-\frac{x}{\tau}}
   \sum_{k=1}^{m+1} (-1)^{m+1+k}\frac{(m+1)!}{k!}
   \binom{m}{k-1}
     \mleft(\frac{x}{\tau}\mright)^k
     \\
     &=&
     \tau^{-m-1} \mye^{-\frac{x}{\tau}}
    E_{m+1}\mleft(\frac{x}{\tau}\mright).
        \end{eqnarray}

\clearpage
 \section*{Appendix B. Proof of equality \eqref{defFsum2}}
 \setcounter{section}{2}
\setcounter{equation}{0}
\numberwithin{equation}{section}
 
 We need to show that \eqref{defmu} is equal to
 \eqref{defFmu} with $F_{d,m}(x)$ defined by \eqref{defFsum2}.
 To this end we prove a more general equality,
 \begin{eqnarray}\label{genidR}  
  \frac{\myd^m}{\myd \tau^m} 
  \left(
  \tau^{
   d +1/2} 
  \mye^{-{x}{\tau}}\right)
  &=&   
  \tau^{ d -m+1/2} \mye^{-{x}{\tau}}
    F_{d,m}\mleft({x}{\tau}\mright).
      \end{eqnarray}

The proof is by induction on $m$. Case $m=0$ is trivial.
  Further  we have:
\begin{eqnarray}
\label{B2}  
 \lefteqn{ \frac{\myd^{m+1}}{\myd \tau^{m+1}}  
  \left( \tau^{
   d +1/2} 
  \mye^{-{x}{\tau}}\right)
  =}
  \\
\label{B3}  
  &=&
   \left(\frac{\myd^m}{\myd \tau^m}  
 \left( \tau^{
   d +1/2} 
  \mye^{-{x}{\tau}}\right)\right)'
  \\
\label{B4}  
 &=&   
     \left(\tau^{ d -m+1/2} \mye^{-{x}{\tau}}
    F_{d,m}\mleft({x}{\tau}\mright)\right)'
\\
 \label{B5}  
    &=&    \left(
      \mye^{-{x}{\tau}}
   \sum_{k=0}^m(-1)^{k}\binom{m}{k}
 \left(d+\tfrac{1}{2}\right)^{\underline{m-k}}
  {x}^k{\tau}^{d-m+k+1/2}\right)'
 \\&=& 
 \label{B6}
  -x
 \mye^{-{x}{\tau}}
   \sum_{k=0}^m(-1)^{k}\binom{m}{k}
 \left(d+\tfrac{1}{2}\right)^{\underline{m-k}}
  {x}^k{\tau}^{d-m+k+1/2}
 \\&&
 \label{B7}
 \quad + 
  \mye^{-{x}{\tau}}
   \sum_{k=0}^m(-1)^{k}({d-m+k+1/2})\binom{m}{k}
 \left(d+\tfrac{1}{2}\right)^{\underline{m-k}}
  {x}^k{\tau}^{d-m+k-1/2}
 \\&=&
 \label{B8}
  -x
 \mye^{-{x}{\tau}}
   \sum_{k=1}^{m+1}(-1)^{k-1}\binom{m}{k-1}
 \left(d+\tfrac{1}{2}\right)^{\underline{m-k+1}}
  {x}^{k-1}{\tau}^{d-m+k-1/2}
 \\&&
 \label{B9}
 \quad + 
  \mye^{-{x}{\tau}}
   \sum_{k=0}^{m+1}(-1)^{k}({d-m+k+1/2})\binom{m}{k}
 \left(d+\tfrac{1}{2}\right)^{\underline{m-k}}
  {x}^k{\tau}^{d-m+k-1/2}
 \\
  \label{B12}
  &=&
  \tau^{ d -m-1/2} \mye^{-{x}{\tau}}
   \sum_{k=0}^{m+1}(-1)^{k}\binom{m+1}{k}
 \left(d+\tfrac{1}{2}\right)^{\underline{m-k+1}}
 \mleft( {x}{\tau}\mright)^k
   \\
       \label{B14}
     &=&
     \tau^{ d -m-1/2} \mye^{-{x}{\tau}}
    F_{d,m+1}\mleft({x}{\tau}\mright).
        \end{eqnarray}

\clearpage
 \section*{Appendix C. Proof of the non-singularity of matrix
 \eqref{defB}}
 
\setcounter{section}{3}
\setcounter{equation}{0}
\numberwithin{equation}{section}
To justify the definition \eqref{Binv} of 
$\alpha_{N,i,j}(d,q)$ and $\beta_{N,i,j}(d,q)$,
we need to prove that the  matrix $B_N(d,q)$ (defined by 
\eqref{defB}) 
is not degenerate. For our purposes we can assume that 
$d$ is a rational number.

Let us consider matrix
\begin{equation}\label{defhatB}
\hat{B}_N(d,x)=\left[
   \begin{array}{ccc}
      \hat{\lambda}_{N,0}(d,x)
           &\dots&
      \hat{\lambda}_{N,2N-1}(d,x)
     \\\vdots&\ddots&\vdots\\
         \hat{\lambda}_{1,0}(d,x)
           &\dots&
      \hat\lambda_{1,2N-1}(d,x)
      \\  \hat\mu_{1,0}(d,x)
           &\dots&
      \hat\mu_{1,2N-1}(d,x)
     \\\vdots&\ddots&\vdots\\
        \hat\mu_{N,0}(d,x)
           &\dots&
      \hat\mu_{N,2N-1}(d,x)
   \end{array}
   \right]
 \end{equation}
 where
$ \hat\lambda_{n,m}(d,x)= E_m({ n  ^2}{x})$
and  $  \hat\mu_{n,m}(d,x)=   F_{d,m}({ n  ^2}{x})$.
 We have:
 \begin{equation}
   \lambda_{n,m}(d,q)=n^d \mye^{-\frac{\ \pi n  ^2}{q}}\hat\lambda_{n,m}(d,\pi/q),
 \end{equation}
 \begin{equation}
   \mu_{n,m}(d,q)=n^d \mye^{-\frac{\ \pi n  ^2}{q}}\hat\mu_{n,m}(d,\pi/q).
 \end{equation}
This implies that
\begin{equation}\label{detdet}
  \det\big({B}_N(d,q)\big)=N!^{2d}\mye^{-\frac{\ \pi  N (N+1) ( 2N+1)}{3q}}\det\big(\hat{B}_N(d,\pi/q)\big)
\end{equation}
and hence it is sufficient to show that the matrix 
$\hat{B}_N(d,\pi/q)$ is not degenerate.

The determinant $\det\big(\hat{B}_N(d,x)\big)$ is a polynomial in $x$.
Both $E_m({{ n  ^2}{x}})$ and $F_{d,m}({{ n  ^2}{x}})$ are polynomials of degree $m$ in $x$,
so the degree of $\det\big(\hat{B}_N(d,x)\big)$ is at most $\sum_{m=0}^{2N-1}m=N(2N-1)$. Therefore,
\begin{equation}\label{rr}
  \det\big(\hat{B}_N(d,x)\big)=\sum_{k=0}^M r_k x^k
\end{equation}
where $M=N(2N-1)$ and $r_0,\dots,r_M$ are rational numbers
(under our assumption that $d$ is rational). 

Let us check that
$r_M\ne 0$.
Calculating this number, we can replace
polynomials  $\hat\lambda_{n,m}(d,x)$ 
and $\hat\mu_{n,m}(d,x)$ in
\eqref{defhatB} with their leading coefficients:
\begin{equation}
  r_M=\det\big(\hat{B}^\mathrm{lc}_N(d)\big)
\end{equation}
where
\begin{equation}
\hat{B}_N^ \mathrm{lc}(d)=\left[
   \begin{array}{ccc}
      \mathrm{lc}\big(\hat{\lambda}_{N,0}(d,x))\big)
           &\dots&
      \mathrm{lc}\big(\hat{\lambda}_{N,2N-1}(d,x)\big)
     \\\vdots&\ddots&\vdots\\
         \mathrm{lc}\big(\hat{\lambda}_{1,0}(d,x)\big)
           &\dots&
      \mathrm{lc}\big(\hat\lambda_{1,2N-1}(d,x)\big)
      \\  \mathrm{lc}\big(\hat\mu_{1,0}(d,x)\big)
           &\dots&
      \mathrm{lc}\big(\hat\mu_{1,2N-1}(d,x)\big)
     \\\vdots&\ddots&\vdots\\
        \mathrm{lc}\big(\hat\mu_{N,0}(d,x)\big)
           &\dots&
      \mathrm{lc}\big(\hat\mu_{N,2N-1}(d,x)\big)
   \end{array},
   \right]
 \end{equation}
 and $ \mathrm{lc}(P)$ is the leading coefficient of the polynomial~$P$.
 According to \eqref{defEsum} and \eqref{defFsum2},
\begin{eqnarray}
 \mathrm{lc}\big(\hat{\lambda}_{n,m}(d,x))\big)&=&\left({n^2}\right)^m, \\
  \mathrm{lc}\big(\hat{\mu}_{n,m}(d,x))\big)&=&\left(-{n^2}\right)^m. 
\end{eqnarray}
Thus, $\hat{B}_N^ \mathrm{lc}(d)$ is a Vandermonde matrix, and hence its determinant, which is equal to $r_M$, is non-zero.

Since $\pi$ is a transcendental number, the polynomial \eqref{rr} 
does not vanish when $x$ is equal to any non-zero rational multiple of~$\pi$.
Therefore, according to \eqref{detdet}, 
matrix $B_N(d,q)$ is non-singular for $q=1,2,\dots$\ .

 \textbf{Remark.}   While polynomials \eqref{rr} 
do not vanish for $x$ equal to a non-zero rational multiple of~$\pi$,
these polynomials seem to have  zeros close to 
$\pi/q$ for $q=4$, $5$, \dots\ . To observe this phenomenon, 
one can calculate the zeros of the polynomials 
\begin{equation}\label{rrr}
  q^M\det\left(\hat{B}_N\left(0,\frac{\pi}{q}\right)\right)=
  \sum_{k=0}^M r_k \pi^k q^{M-k}.
\end{equation}
Table \ref{zerorr}  presents 40 least positive zeros of 
\eqref{rrr} for  $N=12$, $ 15$, and $20$.

\begin{table}[h]\center
   \begin{tabular}{||l|l|l||}
  \hhline{|t:===:t|}
    
  \multicolumn{1}{||c}{$N=12$}
&
  \multicolumn{1}{|c}{$N=15$}
&
  \multicolumn{1}{|c||}{$N=20$}
\\

  \hhline{|:===:|}
    $\  4+1.43487\dotsc\cdot 10^{-39}$
  &
    $\  4-1.40426\dotsc\cdot 10^{-73}$
  &
    $\  4+1.53607\dotsc\cdot 10^{-119}$
  \\ 

    $\  5+3.19394\dotsc\cdot 10^{-27}$
  &
    $\  5+4.91675\dotsc\cdot 10^{-46}$
  &
    $\  5-2.23138\dotsc\cdot 10^{-88}$
  \\ 

    $\  6+5.51749\dotsc\cdot 10^{-22}$
  &
    $\  6+3.25396\dotsc\cdot 10^{-36}$
  &
    $\  6-1.34398\dotsc\cdot 10^{-67}$
  \\ 

    $\  6+4.73456\dotsc\cdot 10^{-19}$
  &
    $\  6+4.89556\dotsc\cdot 10^{-34}$
  &
    $\  6+4.55472\dotsc\cdot 10^{-70}$
  \\ 

    $\  7-2.65631\dotsc\cdot 10^{-13}$
  &
    $\  7-1.39365\dotsc\cdot 10^{-25}$
  &
    $\  7-9.60329\dotsc\cdot 10^{-59}$
  \\ 

    $\  7-1.08677\dotsc\cdot 10^{-15}$
  &
    $\  7-4.45839\dotsc\cdot 10^{-27}$
  &
    $\  7-9.64307\dotsc\cdot 10^{-62}$
  \\ 

    $\  8-2.25893\dotsc\cdot 10^{-9}$
  &
    $\  8-5.06780\dotsc\cdot 10^{-22}$
  &
    $\  8-9.70568\dotsc\cdot 10^{-42}$
  \\ 

    $\  8+8.04754\dotsc\cdot 10^{-16}$
  &
    $\  8-2.23904\dotsc\cdot 10^{-24}$
  &
    $\  8+4.08603\dotsc\cdot 10^{-52}$
  \\ 

    $\  8+3.95665\dotsc\cdot 10^{-11}$
  &
    $\  8-2.18157\dotsc\cdot 10^{-31}$
  &
    $\  8+4.29936\dotsc\cdot 10^{-44}$
  \\ 

    $\  9-1.92047\dotsc\cdot 10^{-7}$
  &
    $\  9-3.48374\dotsc\cdot 10^{-15}$
  &
    $\  9-3.76566\dotsc\cdot 10^{-39}$
  \\ 

    $\  9+1.94667\dotsc\cdot 10^{-13}$
  &
    $\  9-9.75265\dotsc\cdot 10^{-16}$
  &
    $\  9-1.71319\dotsc\cdot 10^{-41}$
  \\ 

    $\  9+2.77801\dotsc\cdot 10^{-6}$
  &
    $\  9+1.16528\dotsc\cdot 10^{-18}$
  &
    $\  9+3.50304\dotsc\cdot 10^{-33}$
  \\ 

    $ 10-3.53263\dotsc\cdot 10^{-4}$
  &
    $ 10-1.66989\dotsc\cdot 10^{-11}$
  &
    $ 10+1.82065\dotsc\cdot 10^{-37}$
  \\ 

    $ 10-1.74943\dotsc\cdot 10^{-5}$
  &
    $ 10-6.70672\dotsc\cdot 10^{-19}$
  &
    $ 10+4.21558\dotsc\cdot 10^{-34}$
  \\ 

    $ 10+8.52601\dotsc\cdot 10^{-11}$
  &
    $ 10+1.65192\dotsc\cdot 10^{-15}$
  &
    $ 10+4.34898\dotsc\cdot 10^{-28}$
  \\ 

    $ 10+8.12171\dotsc\cdot 10^{-8}$
  &
    $ 10+3.96915\dotsc\cdot 10^{-10}$
  &
    $ 10+5.83012\dotsc\cdot 10^{-26}$
  \\ 

    $ 11-5.07370\dotsc\cdot 10^{-3}$
  &
    $ 11-3.58106\dotsc\cdot 10^{-7}$
  &
    $ 11-2.04383\dotsc\cdot 10^{-20}$
  \\ 

    $ 11-3.85253\dotsc\cdot 10^{-6}$
  &
    $ 11-5.39248\dotsc\cdot 10^{-14}$
  &
    $ 11-2.38453\dotsc\cdot 10^{-22}$
  \\ 

    $ 11+4.01883\dotsc\cdot 10^{-7}$
  &
    $ 11+1.43976\dotsc\cdot 10^{-15}$
  &
    $ 11-9.40356\dotsc\cdot 10^{-35}$
  \\ 

    $ 11+6.08232\dotsc\cdot 10^{-3}$
  &
    $ 11+3.10096\dotsc\cdot 10^{-8}$
  &
    $ 11-1.06456\dotsc\cdot 10^{-37}$
  \\ 

    $ 12+5.77332\dotsc\cdot 10^{-9}$
  &
    $ 12-2.84041\dotsc\cdot 10^{-6}$
  &
    $ 12-1.67378\dotsc\cdot 10^{-16}$
  \\ 

    $ 12+1.78233\dotsc\cdot 10^{-5}$
  &
    $ 12-1.02577\dotsc\cdot 10^{-8}$
  &
    $ 12-1.18415\dotsc\cdot 10^{-17}$
  \\ 

    $ 12+9.82536\dotsc\cdot 10^{-4}$
  &
    $ 12+6.41326\dotsc\cdot 10^{-14}$
  &
    $ 12-6.85463\dotsc\cdot 10^{-23}$
  \\ 

    $ 12+4.44090\dotsc\cdot 10^{-2}$
  &
    $ 12+1.16164\dotsc\cdot 10^{-9}$
  &
    $ 12-2.29403\dotsc\cdot 10^{-29}$
  \\ 

    $ 13-6.36804\dotsc\cdot 10^{-2}$
  &
    $ 12+6.51164\dotsc\cdot 10^{-5}$
  &
    $ 12+1.06921\dotsc\cdot 10^{-23}$
  \\ 

    $ 13-1.01354\dotsc\cdot 10^{-3}$
  &
    $ 13-2.08510\dotsc\cdot 10^{-3}$
  &
    $ 13-7.06669\dotsc\cdot 10^{-14}$
  \\ 

    $ 13-3.03801\dotsc\cdot 10^{-5}$
  &
    $ 13-7.52582\dotsc\cdot 10^{-4}$
  &
    $ 13-2.18236\dotsc\cdot 10^{-19}$
  \\ 

    $ 13-2.11913\dotsc\cdot 10^{-8}$
  &
    $ 13+1.53274\dotsc\cdot 10^{-10}$
  &
    $ 13-4.28641\dotsc\cdot 10^{-24}$
  \\ 

    $ 14+2.03601\dotsc\cdot 10^{-10}$
  &
    $ 13+1.58782\dotsc\cdot 10^{-7}$
  &
    $ 13+8.83051\dotsc\cdot 10^{-18}$
  \\ 

    $ 14+1.09077\dotsc\cdot 10^{-7}$
  &
    $ 13+2.55442\dotsc\cdot 10^{-6}$
  &
    $ 13+1.40217\dotsc\cdot 10^{-12}$
  \\ 

    $ 14+1.46130\dotsc\cdot 10^{-4}$
  &
    $ 14-4.09509\dotsc\cdot 10^{-2}$
  &
    $ 14-1.48535\dotsc\cdot 10^{-9}$
  \\ 

    $ 14+3.24088\dotsc\cdot 10^{-3}$
  &
    $ 14+6.79009\dotsc\cdot 10^{-11}$
  &
    $ 14-1.18676\dotsc\cdot 10^{-14}$
  \\ 

    $ 15-1.19707\dotsc\cdot 10^{-6}$
  &
    $ 14+2.49762\dotsc\cdot 10^{-8}$
  &
    $ 14-1.44754\dotsc\cdot 10^{-19}$
  \\ 

    $ 15-6.55257\dotsc\cdot 10^{-9}$
  &
    $ 14+3.03301\dotsc\cdot 10^{-5}$
  &
    $ 14+3.10819\dotsc\cdot 10^{-21}$
  \\ 

    $ 16+5.46932\dotsc\cdot 10^{-12}$
  &
    $ 14+1.04236\dotsc\cdot 10^{-4}$
  &
    $ 14+7.81810\dotsc\cdot 10^{-16}$
  \\ 

    $ 16+6.34726\dotsc\cdot 10^{-4}$
  &
    $ 14+3.40819\dotsc\cdot 10^{-2}$
  &
    $ 14+1.47459\dotsc\cdot 10^{-10}$
  \\ 

    $ 17+5.58344\dotsc\cdot 10^{-4}$
  &
    $ 15+8.00550\dotsc\cdot 10^{-9}$
  &
    $ 15-5.47130\dotsc\cdot 10^{-8}$
  \\ 

    $ 18+2.82566\dotsc\cdot 10^{-4}$
  &
    $ 15+1.77435\dotsc\cdot 10^{-6}$
  &
    $ 15-2.22946\dotsc\cdot 10^{-11}$
  \\ 

    $ 20-2.91020\dotsc\cdot 10^{-1}$
  &
    $ 15+5.24417\dotsc\cdot 10^{-4}$
  &
    $ 15-5.70439\dotsc\cdot 10^{-22}$
  \\ 

    $ 20-4.37113\dotsc\cdot 10^{-2}$
  &
    $ 15+9.02722\dotsc\cdot 10^{-3}$
  &
    $ 15+2.70141\dotsc\cdot 10^{-20}$
  \\ 

  \hhline{|b:===:b|}
\end{tabular}
 
   \caption{Least positive zeros of polynomial 
\eqref{rrr}}
  \label{zerorr}
\end{table}

\clearpage 

\printbibliography

  \end{document}